%% file: Arxiv_JCSE.tex
\documentclass[11pt]{article}

\input{headings_herzet} 
\usepackage{moreverb}
\usepackage[colorlinks,bookmarksopen,bookmarksnumbered,citecolor=red,urlcolor=red]{hyperref}
\usepackage{tikz}
\usepackage{pgfplots}
\usepackage{kfig}
\usepackage{algorithm}
\usepackage{algpseudocode}
\usepackage{url}

\newcommand\BibTeX{{\rmfamily B\kern-.05em \textsc{i\kern-.025em b}\kern-.08em
T\kern-.1667em\lower.7ex\hbox{E}\kern-.125emX}}

\def\volumeyear{2016}
\input{def_macro}	

\begin{document}


\title{Model Reduction from Partial Observations}

\author{C.~Herzet\textsuperscript{(1)}, P.~H\'eas\textsuperscript{(1)}, A.~Dr\'emeau\textsuperscript{(2)}\\
\small{\textsuperscript{(1)} INRIA, Campus de Beaulieu, 35000 Rennes, France}\\
\small{\textsuperscript{(2)} ENSTA Bretagne, 42 rue Fran\c cois Verny, 29200 Brest, France}
}
\date{}

\maketitle

%
%

This paper deals with model-order reduction of parametric partial differential equations (PPDE). More specifically, we consider the problem of finding a good approximation subspace of the solution manifold of the PPDE when only partial information on the latter is available. We assume that two sources of information are available: \textit{i)} a ``rough'' prior knowledge, taking the form of a manifold containing the target solution manifold;  \textit{ii)} partial linear measurements of the solutions of the PPDE (the term partial refers to the fact that observation operator cannot be inverted).  We provide and study several tools to derive good approximation subspaces from these two sources of information. We first identify the best worst-case performance achievable in this setup and propose simple procedures to approximate the corresponding optimal approximation subspace. We then provide, in a simplified setup, a theoretical analysis relating the achievable reduction performance to the choice of the observation operator and the prior knowledge available on the solution manifold.




\section{Introduction}
Our contribution takes place within the context of model reduction for parametric partial differential equations: 
\begin{align}\label{eq:PPDE}
\mathrm{PDE}(\vech,\param)=0, 
\end{align}
where $\vech$ belongs to a Hilbert space $\Sh$ and $\param\in\paramSet$ is a parameter.  When the solution $\vech(\param)$ of \eqref{eq:PPDE} has to be evaluated for many different values $\param\in\paramSet$, the computational effort  
 may become prohibitive. To circumvent this issue, model reduction intends to simplify the resolution of \eqref{eq:PPDE} by (typically) constraining $\vech$ to belong to some low-dimensional subspace $\Sapprox\subset \Sh$. 
 As a matter of fact, the choice of $\Sapprox$ should be made so that each element of the solution manifold $\Ms= \{\vech(\param) \in \Sh: \param\in\paramSet\}$ is well-approximated by some element of $\Sapprox$. 
 Many techniques have been proposed in the literature to identify such subspaces: Taylor \cite{ZAMM:ZAMM19830630105} or Hermite \cite{NME:NME4436} expansions, 
 proper orthogonal decomposition (POD) \cite{9780511622700}, balanced truncation \cite{Antoulas2005Overview}, reduced basis techniques \cite{Quarteroni2011Certified}, etc. 
 
All the methods  mentioned above presuppose some refined knowledge of the solution manifold $\Ms$. 
 For example, it is  typically assumed that the solution manifold $\Mc$ can be finely sampled \cite{ZAMM:ZAMM19830630105,NME:NME4436,9780511622700,Antoulas2005Overview} or, at least, that the set of parameters $\paramSet$ defining $\Ms$ is perfectly known \cite{Quarteroni2011Certified}. Unfortunately, in practice a refined knowledge of $\Ms$ may not always be available. 
  Nevertheless, in many situations one may have access to some partial\footnote{Here, the term ``partial'' refers to the fact that the measurement operator cannot be inverted.} measurements of the elements of $\Ms$. 
 The main question addressed in this paper is therefore as follows: can we benefit from these partial measurements to (complement our prior knowledge and) compute a good approximation subspace for $\Ms$? 

 In order to provide a precise answer to this question, we assume in the rest of this paper that we have the following two ingredients at our disposal:
\begin{itemize}
\item[\textit{i)}] \textit{a prior manifold $\Mps$}, which collects all the knowledge we have ``a priori'' about $\Mc$. 
 The only constraint we impose on $\Mps$ is to be such that
\begin{align}
\Ms\subseteq \Mps, \label{eq:constraints_prior}
\end{align}
\ie the prior manifold must contain all the elements of the target manifold $\Ms$. 
\item[\textit{ii)}]  \textit{a set of partial observations of the elements of $\Ms$}: we assume that we collect, $\forall\vech\in\Ms$,  a set of noiseless linear measurements:
 \begin{align} \label{eq:obseq}
\kbrace{\scap[\Sobasis_j]{\vech} }_{j=1}^\dimSo,
\end{align}
 where $\kbrace{\Sobasis_j }_{j=1}^\dimSo$ is an orthonormal basis (ONB) of some subspace $\So$ and $\scap[\cdot]{\cdot}$ denotes the inner product in $\Sh$. 
\end{itemize}

The prior information typically derives from some physical considerations and/or constraints we may have about the system under study. We give an example of construction of $\Mps$ in Section~\ref{sec:PracticalImplementation}. 
 The nature of the observations available in practice depends on the experimental setup. In this work, we make the assumption that the measurements can be seen as the outputs of a noiseless linear operator.  The noisy setting is not considered hereafter and left for future work.~We note that if $\dimSo=\dimSh$, the observation operator can be inverted to exactly recover the elements of $\Ms$, leading us back to the standard setup. 
 In the sequel, we will thus assume that $\dimSo<\dimSh$, so that some ambiguity on $\vech$ subsists upon the observation of $\kbrace{\scap[\Sobasis_j]{\vech} }_{j=1}^\dimSo$.
 
Although of clear practical interest, only a few contributions have tackled the problem of model reduction from partial measurements.~To the best of our knowledge, the first paper dealing with this question is due to Everson and Sirovich \cite{Everson1995KarhunenLoeve}.~The authors proposed a methodology, dubbed ``Gappy POD'', constructing an approximation subspace when only some elements of each solution $\vech(\param)$ are observed.~This approach has been applied with success to \eg oceanography in \cite{Beckers2003EOF} or fluid mechanics in \cite{BuiThanh2004Aerodynamic,Gunes2006Gappy,Murray2007Application}. Unfortunately, Gappy POD requires some diversity in the observation operator to work properly: as noted in \cite{Gunes2006Gappy}, this method is doomed to produce poor approximation subspaces as soon as some directions of $\Sh$ are never observed.\footnote{In fact, each direction of ``large variation'' of $\Ms$ should be observed in at least one observation.}~This is for example the case when all the elements of the solution manifold are observed through the \textit{same} partial observation operator as in \eqref{eq:obseq}.~In order to circumvent this issue, prior information about $\Ms$ can be included in the reduction process. This approach was recently used in \cite{Maday2014PBDW,Heas2016Reducedorder,Cui2015Datadriven,Cui2014Likelihoodinformed,SpantiniEtAl2015}. 
In \cite{Maday2014PBDW}, the authors suggested to iteratively enrich the approximation subspace by using ``a posteriori'' estimates of some elements of $\Ms$ (the term ``a posteriori'' refers here  to the fact that the estimates stem from the combination of partial observations and some prior knowledge on $\Ms$). In \cite{Heas2016Reducedorder}, the authors of the present work refined this approach in a Bayesian framework: they proposed to include the uncertainty inherent to the a posteriori estimates in the reduction process. An efficient implementation of this Bayesian reduction strategy relying  on adaptive a posteriori sampling   is proposed  in \cite{Cui2015Datadriven}. The works  \cite{Cui2014Likelihoodinformed,SpantiniEtAl2015} are also closely related to this Bayesian perspective of model reduction:  the authors   determine an optimal low-dimensional approximation of the a posteriori distribution for solving efficiently inverse problems governed by PPDEs.  In this paper, we present a rigorous formulation and justification to these approaches in a deterministic framework.~ In particular,  we provide elements of answer to the following questions: 
\begin{itemize}
\item[\textit{i)}]What is the best performance which can be achieved by combining the information provided by the prior manifold $\Mps$ and the collected observations $\{ \kbrace{ {\scap[\Sobasis_j]{\vech}} }_{j=1}^\dimSo \}_{\vech\in\Ms}$?
\item[\textit{ii)}]Can we characterize this ideal performance as a function of the choice of the prior manifold $\Mps$ and the measurement subspace $\So$?  
\item[\textit{iii)}]How to compute a good approximation subspace from $\Mps$ and $\{\kbrace{\scap[\Sobasis_j]{\vech} }_{j=1}^\dimSo\}_{\vech\in\Ms}$ in practice?
\end{itemize}


 As suggested by our above discussion, distinct from, but related to model reduction from partial measurements, is the question of deriving a good estimate of some unknown $\vech(\param)$ by exploiting both a reduced-order model and some collected data.~This paradigm has recently been explored in several papers~\cite{Dihlmann2016Reduced,NME:NME4747,Maday2014PBDW,Binev2015Data}.~In~\cite{Dihlmann2016Reduced}, the authors proposed a reduced version of a Kalman filter and showed that the error on the estimate delivered by the latter can be bounded by a function of the data residual.  In \cite{NME:NME4747} (resp. \cite{Maday2014PBDW}), Maday \etal~considered a data assimilation problem from noiseless (resp.~noisy) linear observations, where the prior model is defined as a low-dimensional approximation subspace of $\Ms$. The same setup was discussed in \cite{Binev2015Data} by Binev \etal  and extended to priors defined as an intersection of degenerate ellipsoids. Some of the theoretical considerations exposed in the present work are built upon the arguments derived in that paper. 
 
The rest of the paper is organized as follows. 
 In Section~\ref{sec:WCOMR} we show that the worst-case optimal performance achievable in our partially-informed setup is characterized by the ``Kolmogorov width'' of some well-defined set.~In Section~\ref{sec:PracticalImplementation}, we propose a simple practical scheme to approximate this worst-case optimal performance and illustrate its performance on the well-known ``thermal-block problem'' in Section \ref{sec:simuresults}.~Finally, in Section \ref{sec:thresults} we provide a theoretical result relating the achievable reduction performance  to the choice of the prior and observation operators. The main steps of the proof of this result are exposed in Section~\ref{sec:reconstruction_guaran} whereas the technical details are postponed to Appendix \ref{app:caractMpost}.  In Appendix \ref{app:wformualtion}, we discuss the weak and algebraic formulations of the PPDE used in our simulations. In order to ease the reading of our work, we also provide a summary of the main notations and expressions useful to understand the paper in Appendix \ref{app:summary}. 


\section{Notational Conventions}

Except if otherwise stated, the  notational conventions used in this paper are as follows.  Italic lowercase boldface letters (as \eg $\vech$, $\Spsbasis$, $\Sobasis$, etc.)  denote elements of the Hilbert space $\Sh$. Uppercase italic letters (\eg $\Sapprox$, $\Sp$, $\Sps$, $\So$, etc.) are used for subspaces of $\Sh$. Lowercase (\eg $\n$, $\x$, etc.) and uppercase boldface (\eg $\X$, $\S$, etc.) letters respectively stand for vector and matrix notations. Lowercase italic  (\eg $a$, $b$, $\alpha$, $\beta$, etc.) and uppercase normal (\eg $\mathrm{N}$, $\mathrm{T}$, etc.) letters  denote scalars. For matrices and vectors, the superscript $\ktranspose{.}$ means transposition. The element located at row $i$ and column $j$ of a matrix $\X$ is denoted $x_{ij}$. 

The inner product and the induced norm associated to $\Sh$ are denoted by $\scap[\cdot]{\cdot}$ and $\kvvbar{ \cdot }$ respectively.~ The space orthogonal to a subspace $X$ with respect to the inner product $\scap[\cdot]{\cdot}$ is written as $X^\perp$.~The subspace induced by a set $\kbrace{\Spsbasis_i}_{i=1}^\dimSps$ is denoted $\spa[\kbrace{\Spsbasis_i}_{i=1}^\dimSps]$.~The distance between an element $\vech\in\Hc$ and a closed subspace $\Sapprox$ is defined as 
\begin{align}
\dist[\vech]{\Sapprox} \triangleq \min_{\vech'\in \Sapprox} \kvvbar{\vech-\vech'}, \nonumber
\end{align}
and 
the projection of $\vech$ onto $\Sapprox$ as 
\begin{align}
\projector[\Sapprox](\vech) \triangleq \kargmin_{\vech'\in \Sapprox}\kvvbar{\vech-\vech'}. \nonumber
\end{align}
The notation  $\projector[\Sapprox](\Sps)$ stands for the set $\kbrace{\projector[\Sapprox](\vech): \vech\in\Sps}$. The operator $\oplus$ is used to denote the direct sum between two subsets of $\Sh$, \eg  $\Sapprox \oplus \Sps= \kbrace{{\boldsymbol s} + \Spsbasis: {\boldsymbol s}\in \Sapprox, \Spsbasis\in \Sps}$.~Finally, $B_{\Spwidth}=\kbrace{\vech\in\Sh: \kvvbar{\vech}\leq\Spwidth}$ is the $\kvvbar{\cdot}$-ball of radius $\Spwidth$.

\section{Worst-case Optimal Model Reduction}\label{sec:WCOMR}

In this section, we tackle the problem of finding a good approximation subspace from a worst-case perspective.~In Section \ref{sec:defconstraints}, we first discuss the nature of the information provided by the prior manifold \eqref{eq:constraints_prior} and the collected observations \eqref{eq:obseq}, and identify the set of manifolds compatible with the latter.~In Section \ref{sec:WCROM}, we then characterize the best worst-case performance achievable in our partially-informed setup.


\subsection{Feasible and Posterior Manifolds}\label{sec:defconstraints}

The prior manifold and the partial measurements 
 provide some valuable information about the unknown manifold $\Ms$ since they both define a set of constraints (discussed below) which are known to be satisfied by the latter.~However, $\Ms$ is usually not the only manifold verifying these constraints. 
 In the sequel, we will denote the set of manifolds satisfying the constraints imposed by  $\Mps$ and the partial measurements $\{\kbrace{\scap[\Sobasis_j]{\vech}  }_{j=1}^\dimSo\}_{\vech\in\Ms}$ as $\feasmanifoldset$; the elements of $\feasmanifoldset$ will be referred to as ``feasible'' manifolds.

  In the rest of this subsection, we give a precise characterization of $\feasmanifoldset$ and emphasize that the largest\footnote{The term ``largest'' must be understood as follows: if $\Mstrial\in \feasmanifoldset$ then $\Mstrial \subseteq \Mpost$.} element of this set is defined by 
\begin{align}\label{eq:defMpost}
\Mpost = \Mps \cap \kparen{\cup_{\vech\in\Ms} \Hplan_\vech},
\end{align}
where 
\begin{align}
\Hplan_\vech & \triangleq \kbrace{\vech'= \vech +  \Sobasis^\perp: \Sobasis^\perp\in \So^\perp}.\nonumber
\end{align}
This observation will turn out to be crucial in the next subsection to characterize the worst-case optimal approximation subspace. 

Let us first discuss the constraints defining $\feasmanifoldset$.~First, it is clear that any manifold $\Mstrial$ compatible with our prior assumption \eqref{eq:constraints_prior} should be such that
\begin{align}\label{eq:PriorConstraint}
\Mstrial \subseteq \Mps. 
\end{align}
Secondly,  a manifold $\Mstrial$ compatible with the received observations should  reproduce exactly the same set of measurements as those obtained from $\Ms$ when measured with the same observation operator. More specifically, any manifold $\Mstrial$ compatible with the received observations should be such that
\begin{align}\label{eq:MeasureConstraint}
\kbrace{\kbrace{\scap[\Sobasis_j]{\vech} }_{j=1}^\dimSo}_{\vech\in\Mstrial}=\kbrace{\kbrace{\scap[\Sobasis_j]{\vech} }_{j=1}^\dimSo}_{\vech\in\Ms}.
\end{align}
We thus define the set of feasible manifolds as 
\begin{align}\label{eq:defFeasibleSet}
\feasmanifoldset = \kbrace{\Mstrial : \mbox{\eqref{eq:PriorConstraint} and \eqref{eq:MeasureConstraint} hold}}.
\end{align}
In other words, $\feasmanifoldset$ represents the set of manifolds which are compatible with both our prior assumption \eqref{eq:constraints_prior} and the set of collected observations. It is obvious from \eqref{eq:PriorConstraint}-\eqref{eq:defFeasibleSet} that $\Ms\in\feasmanifoldset$.

We now prove the following lemma:
\begin{lemma}\label{lemma:Mpost}
$\Mpost$ is the largest element of $\feasmanifoldset$.
\end{lemma}
\textit{Proof:} We first show that
\begin{align}\label{eq:necessarycondfeas}
\left\{
\begin{array}{l}
\Mstrial \subseteq \Mps,\\ 
\Mstrial \subseteq \cup_{\vech\in\Ms} \Hplan_\vech,
\end{array}
\right.
\end{align}
are necessary conditions for $\Mstrial\in\feasmanifoldset$ and, \eqref{eq:necessarycondfeas} together with 
\begin{align}\label{eq:sufficientcond}
\Ms\subseteq\Mstrial,
\end{align}
is a sufficient condition for $\Mstrial\in\feasmanifoldset$. The proof of the main result will derive straightforwardly from these conditions.  

The necessity of $\Mstrial \subseteq \Mps$ is obvious from the definition of $\feasmanifoldset$. The necessity of $\Mstrial \subseteq \cup_{\vech\in\Ms} \Hplan_\vech$ can be shown as follows. First, note that the set of elements of $\Sh$ leading to a given set of observations  $\{\scap[\Sobasis_j]{\vech}\}_{j=1}^\dimSo$ is an affine subspace  defined as
\begin{align}
\Hplan_\vech 
&= \kbrace{\vech' : \langle \Sobasis_j, \vech' \rangle =\langle \Sobasis_j, \vech \rangle \mbox{ for $j=1,\ldots, m$}}, \nonumber\\
&= \kbrace{\vech'= \vech +  \Sobasis^\perp: \Sobasis^\perp\in \So^\perp}. \nonumber
\end{align}
Hence, if $\Mstrial \nsubseteq \cup_{\vech\in\Ms} \Hplan_\vech$, then $\kexists{\vechtrial\in\Mstrial}$ such that
\begin{align}
\kbrace{\scap[\Sobasis_j]{\vechtrial} }_{j=1}^\dimSo \notin \kbrace{\kbrace{\scap[\Sobasis_j]{\vech} }_{j=1}^\dimSo}_{\vech\in\Ms}\ , \nonumber
\end{align}
and therefore \eqref{eq:MeasureConstraint} cannot be satisfied. The necessity of $\Mstrial \subseteq \cup_{\vech\in\Ms} \Hplan_\vech$  is thus obtained  by contraposition.

We now show that \eqref{eq:necessarycondfeas}-\eqref{eq:sufficientcond} is sufficient for $\Mstrial\in\feasmanifoldset$. Since the first condition in \eqref{eq:necessarycondfeas} is identical to \eqref{eq:PriorConstraint}, we only need to show that \eqref{eq:necessarycondfeas}-\eqref{eq:sufficientcond} implies \eqref{eq:MeasureConstraint}.~First, from our previous discussion, we have that $\Mstrial \subseteq \cup_{\vech\in\Ms} \Hplan_\vech$ implies that 
\begin{align}
\kbrace{\kbrace{\scap[\Sobasis_j]{\vech} }_{j=1}^\dimSo}_{\vech\in\Mstrial}\subseteq\kbrace{\kbrace{\scap[\Sobasis_j]{\vech} }_{j=1}^\dimSo}_{\vech\in\Ms}.\nonumber
\end{align}
Moreover, if $\Ms\subseteq\Mstrial$ holds then we also have
\begin{align}
\kbrace{\kbrace{\scap[\Sobasis_j]{\vech} }_{j=1}^\dimSo}_{\vech\in\Ms}\subseteq\kbrace{\kbrace{\scap[\Sobasis_j]{\vech} }_{j=1}^\dimSo}_{\vech\in\Mstrial}.\nonumber
\end{align}
Combining the last two inclusions, we obtain \eqref{eq:MeasureConstraint}. 

We finally prove the statement of Lemma \ref{lemma:Mpost} by exploiting the necessary 
 and sufficient conditions 
  defined above. Since any $\Mstrial\in\feasmanifoldset$ must satisfy \eqref{eq:necessarycondfeas}, we have
\begin{align}
\Mstrial \subseteq \Mps \cap \kparen{\cup_{\vech\in\Ms} \Hplan_\vech} = \Mpost. \nonumber
\end{align}
It thus remains to show that $\Mpost\in \feasmanifoldset$. Now, $\Mpost$ verifies \eqref{eq:necessarycondfeas} by definition.~Moreover, \eqref{eq:sufficientcond} also holds because $\Ms$ is included in both $\Mps$ and $\cup_{\vech\in\Ms} \Hplan_\vech$. Hence, $\Mpost\in \feasmanifoldset$ by virtue of the sufficiency of \eqref{eq:necessarycondfeas}-\eqref{eq:sufficientcond}.  
\hfill $\square$\\

$\Mpost$ summarizes the uncertainty about the unknown manifold $\Ms$ by gathering the information provided by the prior model and the partial measurements.~In particular, it is the smallest subset of $\Sh$ containing all the manifolds compatible with the prior constraint \eqref{eq:constraints_prior} and the received observations. 
In the sequel, we will thus refer to $\Mpost$ as the ``posterior'' manifold because of its analogy with the posterior probabilities defined in a Bayesian framework: both characterize the uncertainty remaining on some quantity of interest upon the combination of some prior and observation models \cite{Robert2001Bayesian}. 
%

\subsection{Worst-case Optimal Model Reduction}\label{sec:WCROM}
In the model-reduction literature, an ideal figure of merit to assess the reducibility of a manifold $\Ms$ is its {``Kolmogorov $i$-width''}: 
\begin{align}\label{eq:KolmogorovWith}
\Kspec_i(\Ms) &= \inf_{\Sapprox:\dim(\Sapprox)=i} \kparen{\sup_{\vech\in\Ms} \dist[\vech]{S}}. 
\end{align}
It is clear from its definition that  $\Kspec_i(\Ms)$ provides the best  worst-case  error achievable by any approximation subspace of dimension $i$. If the infimum of \eqref{eq:KolmogorovWith} can be attained\footnote{In order to ease the discussion, we will always assume hereafter that all the suprema and infima exist and can be attained. 
 We refer the reader to \cite{Pinkus1985NWidths} for a detailed discussion on the conditions ensuring  the existence of the extremizers of the Kolmogorov $i$-width.}, a worst-case optimal approximation subspace is thus given by
\begin{align}\label{eq:FIWO}
\Sperf_i &\in \kargmin_{\Sapprox:\dim(\Sapprox)=i} \kparen{\sup_{\vech\in\Ms} \dist[\vech]{S}}.
\end{align}
The resolution of this problem obviously entails the knowledge of the manifold to reduce, \ie $\Mc$.~In the setup considered in this paper, which presupposes that $\Ms$ is unknown,  computing an approximation subspace according to \eqref{eq:FIWO} is therefore not possible. 

Nevertheless, as discussed in the previous subsection, the presence of prior information and partial measurements on the unknown manifold $\Ms$ reduces the uncertainty about its localization in $\Sh$. More specifically, any manifold compatible with the prior information and the received observations must belong to the feasible set $\feasmanifoldset$.~A sensible approach, followed hereafter, then consists of including this information within the worst-case criterion used to evaluate the approximation subspace. More precisely, we suggest to compute the approximation subspace as the solution of
\begin{align}\label{eq:WCoptimization1}
\Spost_i&\in \kargmin_{\Sapprox : \dim(\Sapprox)=i} \kparen{\sup_{\Mstrial\in\feasmanifoldset} \sup_{\vech\in\Mstrial} \dist[\vech]{\Sapprox}}. 
\end{align}
Problem \eqref{eq:WCoptimization1} is tantamount to finding the approximation subspace minimizing the maximum approximation error over all the feasible manifolds.~We note that since the posterior manifold $\Mpost$ defined in \eqref{eq:defMpost} is the largest element of $\feasmanifoldset$ (see Lemma \ref{lemma:Mpost}), \eqref{eq:WCoptimization1} can also be rewritten as
\begin{align}
\Spost_i &\in \kargmin_{\Sapprox: \dim(\Sapprox)=i}\kparen{\sup_{\vech\in\Mpost} \dist[\vech]{\Sapprox}}.\label{eq:MainProblem}
\end{align}
By definition of $\Spost_i$ and since $\Ms\subseteq\Mpost\subseteq \Mps$, we have 
\begin{align}
\Kspec_i(\Ms) \leq \sup_{\vech\in\Ms} \dist[\vech]{\Spost_i} \leq \Kspec_i(\Mpost)\leq \Kspec_i(\Mps).  \label{eq:upperBKv} 
\end{align}
The first inequality follows from the definition the Kolmogorov $i$-width of $\Ms$, the second from $\Ms\subseteq\Mpost$ and the last one is a consequence of $\Mpost\subseteq \Mps$. 

In the light of \eqref{eq:upperBKv}, we see that in the partially-informed setup considered in this paper, the optimal reduction performance is lower bounded by $\Kspec_i(\Ms)$ and upper bounded by $\Kspec_i(\Mpost)$. 
 The gap between $\Kspec_i(\Ms)$ and $\Kspec_i(\Mpost)$ ``materializes'' the loss of reducibility which can occur by working in a partially-informed setting rather than a perfectly-informed one. 
 We also note that $\Kspec_i(\Mps)$ characterizes the best achievable worst-case performance when the prior constraint \eqref{eq:constraints_prior} is the only information available to the practitioner (\ie there are no observations). More precisely we have
\begin{align}
\inf_{\Sapprox : \dim(\Sapprox)=i} \kparen{\sup_{\Mstrial\subseteq\Mps} \sup_{\vech\in\Mstrial} \dist[\vech]{\Sapprox}}
= \inf_{\Sapprox : \dim(\Sapprox)=i} \kparen{\sup_{\vech\in\Mps} \dist[\vech]{\Sapprox}}=\Kspec_i\kparen{\Mps}. \nonumber
\end{align}
As expected, from a worst-case perspective, there is thus always a gain in exploiting the received observations on top of the prior information; the gain brought by the former is characterized by the gap between $\Kspec_i\kparen{\Mpost}$ and $\Kspec_i\kparen{\Mps}$. In Section~\ref{sec:thresults}, we will discuss more specifically the connections between $\Kspec_i\kparen{\Ms}$, $\Kspec_i\kparen{\Mpost}$ and $\Kspec_i\kparen{\Mps}$ in a simplified setup. 

\begin{algorithm}[t]
\begin{algorithmic}[0]
\State \textbf{inputs}: $\kbrace{\vech_j}_{j=1}^\nsnapshots$
\State \textbf{init}: $\hat{\Sapprox}=\{\0\}$
\While{Stopping criterion is not satisfied}
	\State Compute $\vech = \kargmax_{\vechtrial\in \kbrace{\vech_j}_{j=1}^\nsnapshots}\kvvbar{\vechtrial-\projector[\hat{\Sapprox}](\vechtrial)}$
	\State Set 	$\Sapproxa = \Sapproxa \oplus \spa[\vech]$
\EndWhile
\State \textbf{ouput}: $\Sapproxa$
\end{algorithmic}
\caption{Greedy algorithm\label{alg:Greedy}}
\end{algorithm}

\section{Practical Implementation: Greedy Procedure and Sampling Schemes}\label{sec:PracticalImplementation}

Solving \eqref{eq:MainProblem} is typically an intractable problem. Suboptimal approaches have therefore to be considered to find good approximated solutions.~Interestingly, since \eqref{eq:MainProblem} shares exactly the same structure as \eqref{eq:FIWO} (with the difference that the supremum is taken over $\Mpost$ rather than $\Ms$),  one can benefit from the numerous suboptimal techniques which have been proposed in the standard setup to tackle our partially-informed problem.~In the sequel, we will  focus more specifically on the standard ``greedy'' procedure
 described  in Algorithm~\ref{alg:Greedy}, see \eg \cite[Section 7.1.1]{Quarteroni2016Reduced}. The stopping criterion mentioned in the procedure may be for example the dimension of the approximation subspace or some accuracy requirements. 

The algorithm presupposes that a set of elements of the manifold to reduce (commonly referred to as ``snapshots'' in the literature), say $\kbrace{\vech_j}_{j=1}^\nsnapshots$, is available. At each iteration, the procedure increases the dimension of the approximation subspace, by including  one of the snapshots leading to the largest projection error.~The performance of the greedy algorithm thus depends on the choice of $\kbrace{\vech_j}_{j=1}^\nsnapshots$.~A sensible choice, that we will follow hereafter, consists of drawing the snapshots randomly from some distribution supported on $\Mpost$.~We note that
\begin{align}
\Mpost 
&= \Mps \cap \kparen{\cup_{\vech\in\Ms}\Hplan_\vech}\nonumber,\\
&= \cup_{\vech\in\Ms}\kparen{\Mps \cap \Hplan_\vech}\nonumber,
\end{align}
and therefore a sampling of $\Mpost$ can be achieved by drawing snapshots from some distribution supported on $\Mps \cap \Hplan_\vech$, $\forall \vech\in \Ms$. 

The structure of $\Mps \cap \Hplan_\vech$ depends obviously on the choice of $\Mps$. From a practical point of view, some choices of $\Mps$ may thus enable  more favorable implementations than others. In Section~\ref{sec:favorable_priors}, we introduce some specific choices for $\Mps$ which lead to a simple mathematical characterization of $\Mps \cap \Hplan_\vech$. We emphasize that, in many situations, this type of priors  can be obtained quite easily by applying standard reduced-order model techniques.  In Section~\ref{sec:favorable_bases}, we present some material and notations needed to properly characterize $\Mps \cap \Hplan_\vech$. In Section~\ref{eq:samplingStrategies}, we emphasize that there exist simple schemes to sample $\Mpost$ for the specific choices of $\Mps$ introduced in Section \ref{sec:favorable_priors}. Our approach relies on the derivations of Binev \etal showing that $\Mps \cap \Hplan_\vech$ corresponds to a high-dimensional ellipsoid with orthogonal principal axes for some choice of  $\Mps$, see \cite[section 2.4]{Binev2015Data}. In Section~\ref{sec:point-estimate}, we establish some connections between the proposed procedure and some other approaches based on point estimates of the elements of $\Ms$. 
 Finally, in Section~\ref{sec:complexity}, we analyze the complexity of the proposed methods and discuss some implementation issues.
\subsection{Some Specific Choices for $\Mps$} \label{sec:favorable_priors}

In this section, we advocate that the following choice of prior manifold can be obtained from standard model-order reduction techniques:
 \begin{align} \label{eq:PriorUoS0}
\Mps=\cap_{j=1}^{\nsubspace} \kbrace{\vech : \dist[\vech]{\Sps_j}\leq \Spswidth_j}, 
\end{align}
 where $\Sps_j$ is some $\dimSps_j$-dimensional subspace, $\Spswidth_j$ is some positive scalar and $\nsubspace$ is an integer. 
  In Section~\ref{eq:samplingStrategies}, we will emphasize that such a choice for $\Mps$ allows for an easy implementation of the sampling of $\Mps \cap \Hplan_\vech$.

In order to show that \eqref{eq:PriorUoS0} can (for example) be obtained from standard model-order reduction techniques, let us first  consider the case where  the uncertainty one has on the set of possible solutions of the PPDE (\ie $\Ms$) is due to an imperfect knowledge of the set of feasible parameters $\paramSet$ (The general case will be discussed at the end of this section). Although $\paramSet$ may not be precisely known, an information the practitioner usually has at its disposal is that $\paramSet$ is contained in some larger set  $\paramSetp$, \ie $\paramSet\subseteq\paramSetp$. 
 This relaxed set defines another manifold
\begin{align}\label{eq:exMprior}
\Mrs &\triangleq \kbrace{\vech(\param) \in\Sh:\param\in\paramSetp },
\end{align}
which obeys $\Ms \subseteq \Mrs$ (because $\paramSet\subseteq\paramSetp$). 
 Now, constructing a reduced-order model for $\Mrs$ defined as in \eqref{eq:exMprior}  via \eg reduced-basis techniques \cite{Quarteroni2011Certified} leads to a sequence of subspaces 
\begin{align}\label{eq:subspaces}
\Sps_1 \subset \ldots \subset \Sps_\nsubspace,
\end{align}
and scalars
\begin{align}\label{eq:scalars}
\Spswidth_1\geq \ldots \geq \Spswidth_\nsubspace,
\end{align}
such that\footnote{For the sake of precision, let us mention that, as far as reduced-basis techniques are concerned, the inclusion in \eqref{eq:constraintVeps} is guaranteed only if the reduced-basis procedures exploit the \textit{whole} set $\paramSetp$ in the computation of $\Sps_j$ and $\Spswidth_j$. In practice, however, these algorithms only use a fine discretization of $\paramSetp$ to speed up the computations. In such a case, ensuring that  \eqref{eq:constraintVeps} holds requires some additional care. In the sequel, we do not elaborate on this technical issue and assume that $\Sps_j$ and $\Spswidth_j$ are such that the inclusion \eqref{eq:constraintVeps} holds. }
\begin{align}\label{eq:constraintVeps}
{\Mrs \subseteq} \kbrace{\vech : \dist[\vech]{\Sps_j}\leq \Spswidth_j}\ \mbox{ for all $j\leq \nsubspace$}. 
\end{align}
Since $\Ms\subseteq \Mrs$
, one can thus construct a prior of the form \eqref{eq:PriorUoS0}  satisfying \eqref{eq:constraints_prior}, from any combination of the subspaces \eqref{eq:subspaces} and scalars \eqref{eq:scalars}. 

We note that this procedure applies even if $\Mrs$ is not defined via a relaxed set of parameters $\paramSetp$ as in \eqref{eq:exMprior}. In a general setting, the only constraints to satisfy in order to fulfil the condition $\Ms \subseteq \cap_{j=1}^{\nsubspace} \kbrace{\vech : \dist[\vech]{\Sps_j}\leq \Spswidth_j}$ are as follows: \textit{i)} identify a manifold $\Mrs$ such that $\Ms\subseteq\Mrs$; \textit{ii)} apply a model-order reduction technique on $\Mrs$ which certifies that the sequence of subspaces $\Sps_j$ and scalars $\Spswidth_j$ obey \eqref{eq:constraintVeps}.

\subsection{Definition of Suitable Representation Bases} \label{sec:favorable_bases}


In this section, we introduce some representation bases which will be useful for the characterization of $\Mps \cap \Hplan_\vech$ in the next section. 
We let $\So$ and $\Sps$ denote respectively $\dimSo$- and $\dimSps$-dimensional subspaces of $\Hc$. In the next section, $\So$ and $\Sps$ will  play the role of the observation subspace defined in \eqref{eq:obseq} and one of the ``prior'' subspaces $\Sps_j$ appearing in the definition of $\Mps$ in \eqref{eq:PriorUoS0}. 

Let $\G\in\Rbb^{\dimSo\times\dimSps}$ be the matrix representation of the projector $\projector[\So]$ from $\Sps$ to $\So$ in some (arbitrary) ONBs $\kbrace{\Sobasis_j}_{j=1}^\dimSo$ and $\kbrace{\Spsbasis_j}_{j=1}^\dimSps$, \ie 
\begin{align}\nonumber
g_{ij} \triangleq \scap[\Sobasis_i]{\Spsbasis_j}.
\end{align}
 We define new ONBs, $\kbrace{\Sobasis_j^*}_{j=1}^\dimSo$ for $\So$ and $\kbrace{\Spsbasis_j^*}_{j=1}^\dimSps$ for $\Sps$, as
\begin{align}
\Sobasis^*_j& \triangleq \sum_{i=1}^m \Sobasis_i\, x_{ij},\nonumber\\
\Spsbasis^*_j& \triangleq \sum_{i=1}^n \Spsbasis_i\, z_{ij},\nonumber
\end{align}
where $x_{ij}$ and $z_{ij}$ are the coefficients of the matrices appearing in the singular value decomposition of $\G=\X \S \ktranspose{\Z}$; $\X\in\Rbb^{\dimSo\times\dimSo}$ and $\Z\in\Rbb^{\dimSps\times\dimSps}$ are orthogonal matrices; $\S\in\Rbb^{\dimSo\times\dimSps}$ is a rectangular,  diagonal matrix.  

Clearly, from the definition of $\kbrace{\Sobasis_i^*}_{i=1}^\dimSo$ and $\kbrace{\Spsbasis_i^*}_{i=1}^\dimSps$, we have
\begin{align}\label{eq:correlationoptbases}
\scap[\Sobasis^*_i]{\Spsbasis^*_j}&= s_{ij}. 
\end{align}
In the sequel, we will use the shorthand notation $\svdecr_j$ to refer to the $j$th diagonal element of  $\S$, that is $\kbrace{\svdecr_j}_{j=1}^\mindim$ represents the set of singular values of $\G$.\footnote{In order to provide some geometrical interpretation to the reader, let us mention that the singular value $\svdecr_j$ can also be understood as the cosine of the ``angle'' between the singular vectors $\Sobasis^*_j$ and $\Spsbasis^*_j$.} The singular values are assumed to be sorted in a decreasing order of magnitude, \ie
\begin{align}
1\geq \svdecr_{1} \geq \cdots \geq \svdecr_{\mindim} \geq 0.  \nonumber
\end{align}
The first inequality follows from the fact that $\G$ is the matrix representation of a projection operator: we thus necessarily have that $\svdecr_j\leq 1$. Moreover, if $\svdecr_j=1$, we must have $\Sobasis_j^*=\Spsbasis_j^*$. 

We define the following short-hand notations that will be useful in the rest of the paper:
\begin{align}
\begin{array}{ll}
\maxosv & \triangleq  \card[ \kbrace{j :\svdecr_j=1}],\nonumber\\
\maxzsv & \triangleq  \card[ \kbrace{j :\svdecr_j>0}].\nonumber\\
\end{array}
\end{align}
 From an operational point of view, $\maxosv$ represents the number of dimensions of $\Sps$ which are included in $\So$, that is $\maxosv=\dim\kparen{\So\cap\Sps}$.~Moreover, $\dimSps -\maxzsv$ corresponds to the number of dimensions of $\Sps$ which are orthogonal to $\So$, that is $\dimSps-\maxzsv=\dim\kparen{\Sop\cap\Sps}$. In a nutshell, $\maxzsv$ thus represents the number of measurements (out of $\dimSo$) providing information about the position of the points in $\Sps$.

 Finally let us mention that $\Sop$ can be decomposed as the following direct sum of two orthogonal subspaces (see Appendix \ref{sec:charac_MpcapHplan}):
\begin{align}\label{eq:decompW}
\Sop &= \projector[\Sop](\Sps)\, \oplus \Sop\cap\Spsp. 
\end{align}
Moreover an ONB for $\projector[\Sop](\Sps)$ can be expressed in terms of the elements of $\kbrace{\Sobasis_j^*}_{j=1}^\dimSo$ and $\kbrace{\Spsbasis_j^*}_{j=1}^\dimSps$. More specifically, letting 
\begin{align}
\intbasis_j\triangleq \kparen{1-\svdecr_j^2}^{-\frac{1}{2}}\kparen{\Spsbasis^*_j - \svdecr_j \Sobasis^*_j}, \nonumber
\end{align}
we have that  $\kbrace{\intbasis_j}_{j=\maxosv+1}^\maxzsv\cup\kbrace{\Spsbasis_j^*}_{j=\maxzsv+1}^\dimSps$
 forms an ONB of $\projector[\Sop](\Sps)$, see Appendix \ref{sec:charac_MpcapHplan}.
 

\subsection{Sampling Strategies}\label{eq:samplingStrategies}


\begin{algorithm}[t]
\begin{algorithmic}[0]
\State \textbf{inputs}: $\kbrace{\scap[\Sobasis_j]{\vech} }_{j=1}^\dimSo$, $\Spswidth$, $\Sps$\vspace{0.1cm}
\State \textbf{init}: Compute $\kbrace{\Spsbasis_j^*}_{j=1}^\dimSps$, $\kbrace{\Sobasis_j^*}_{j=1}^\dimSo$, $\kbrace{\svdecr_j}_{j=1}^\mindim$ as in Section \ref{sec:favorable_bases}
\State \ \ \ \ \ \ \ \ Evaluate $\kangle{\Sobasis_j^*,\vech} = \sum_{i=1}^\dimSo \scap[\Sobasis_i]{\vech} x_{ij}$ for $j=1,\ldots,\dimSo$ \vspace{0.1cm}
\While{$i< \nsamples$}
	\State 1)  Draw $\pi$ according to some distribution supported on $\kbracket{0,1}$ 
	\State 2) Draw $\gamma$ uniformly on $\kbracket{0,\Spswidth^2-\sum_{j=\maxzsv+1}^\dimSo \kangle{\Sobasis_j^*,\vech}^2}$ 
	\State 3) Draw $\kbrace{b_j}_{\maxosv+1}^{\maxzsv}$ uniformly on the $(\maxzsv-\maxosv)$-dimensional unit ball;
	\State \ \ \ \,  scale the result so that $\sum_{j=\maxosv+1}^{\maxzsv} b_j^2=\gamma\pi$
	\State 4) Draw $\kbrace{d_j}_{\maxzsv+1}^{\dimSps}$ uniformly on  $\Rbb^{\dimSps-\maxzsv}$ 
	\State 5) Draw $\vechVpWp$ uniformly on $\Sps^\perp\cap \So^\perp$; scale the result so that $\kvvbar{\vechVpWp}^2=\gamma\kparen{1-\pi^2}$
	\State 6) Set $\vech_i = \cel - \sum_{j=\maxosv+1}^{\maxzsv} b_j \svdecr_j^{-1}\intbasis_j
+ \sum_{j=\maxzsv+1}^{\dimSps} d_j \Spsbasis^*_j + \vechVpWp$ 
	\State $i=i+1$
\EndWhile \vspace{0.1cm}
\State \textbf{ouput}: $\kbrace{\vech_i \in \Mps \cap \Hplan_\vech}_{i=1}^{\nsamples}$ 
\end{algorithmic}
\caption{Random Sampling of $\Mps \cap \Hplan_\vech$ for $\nsubspace=1$\label{alg:Sampling1}}
\end{algorithm}

We now expose our  strategies to sample $\Mps \cap \Hplan_\vech$ with $\Mps$ defined as in \eqref{eq:PriorUoS0}. 
 First note that each set $\kbrace{\vech : \dist[\vech]{\Sps_j}\leq \Spswidth_j}$ can be seen as a degenerate ellipsoid (the degenerate directions are defined by  the subspace $\Sps_j$);  the definition in \eqref{eq:PriorUoS0} thus corresponds to the intersection of $\nsubspace$ degenerate ellipsoids. In what follows, we will see that this particular geometrical structure allows for an easy sampling of $\Mps \cap \Hplan_\vech$.

We first consider \eqref{eq:PriorUoS0} in the case where $\nsubspace=1$, that is $\Mps$ is defined as
 \begin{align} \label{eq:PriorSS}
\Mps=\kbrace{\vech : \dist[\vech]{\Sps}\leq \Spswidth}, 
\end{align}
for some $\dimSps$-dimensional subspace $\Sps$ and scalar $\Spswidth>0$.~In this setup, $\Mps \cap \Hplan_\vech$ takes the form of a high-dimensional ellipsoid with orthogonal principal axes.~More specifically, Binev \etal \cite[Section 2.4]{Binev2015Data} (see also Appendix \ref{app:def_of_MpscapHplan}) showed that 
 $\Mps \cap \Hplan_\vech$ can be characterized as follows: 
 \begin{align}\label{eq:defMpcapHplan}
\Mps \cap \Hplan_\vech &= \cel \oplus \manifoldperp_\vech,
\end{align}
where 
\begin{align}\label{eq:cel}
\cel &= \sum_{j=1}^\maxzsv \scap[\Sobasis^*_j]{\vech} \svdecr_j^{-1} \Spsbasis^*_j + \sum_{j=\maxzsv+1}^\dimSo \scap[\Sobasis^*_j]{\vech} \Sobasis_j^*,
\end{align}
and $\manifoldperp_\vech$ is an ellipsoid defined as
\begin{align}\label{eq:constrainttot0}
\manifoldperp_\vech &= \left\{
\begin{array}{ll}
 \vechtrial  &=  - \sum_{j=\maxosv+1}^{\maxzsv} b_j \svdecr_j^{-1}\intbasis_j
+ \sum_{j=\maxzsv+1}^{\dimSps} d_j \Spsbasis^*_j + \vechVpWp\\[0.2cm]
 \mbox{with} & \left\{
\begin{array}{l}
\vechVpWp \in \Sps^\perp\cap\So^\perp \\
 \sum_{j= \maxosv+1}^\maxzsv b_j^2 + \kvvbar{\vechVpWp}^2 \leq \Spswidth^2 - \sum_{j=\maxzsv+1}^\dimSo \scap[\Sobasis^*_j]{\vech}^2
\end{array}
\right.
\end{array}
\right\}.
\end{align}


In \eqref{eq:cel}-\eqref{eq:constrainttot0}, we have used the definitions of  $\Spsbasis^*_j$, $\Sobasis^*_j$ and $\intbasis_j$ introduced in Section \ref{sec:favorable_bases} with the following conventions: $\Sps$ is the subspace characterizing $\Mps$ in \eqref{eq:PriorSS}, $\So$ is the observation subspace introduced in \eqref{eq:obseq}. 

From \eqref{eq:defMpcapHplan}-\eqref{eq:constrainttot0} it is clear that   $\Mps \cap \Hplan_\vech$ is an ellipsoid centered in $\cel$.  
The set $\manifoldperp_\vech$ characterizes the deviation of the elements of $\Mps \cap \Hplan_\vech$ from its center.~From our final remark in Section~\ref{sec:favorable_bases} (namely, $\Sop$ can be decomposed as in \eqref{eq:decompW} and $\kbrace{\intbasis_j}_{j=\maxosv+1}^\maxzsv\cup\kbrace{\Spsbasis_j^*}_{j=\maxzsv+1}^\dimSps$ forms an ONB of $\projector[\Sop](\Sps)$), we can deduce that: \textit{i)}  $\manifoldperp_\vech\subseteq \Sop$;  \textit{ii)} the elements $\kbrace{\intbasis_j}_{j=\maxosv+1}^\maxzsv\cup\kbrace{\Spsbasis_j^*}_{j=\maxzsv+1}^\dimSps$ correspond to principal axes of the ellipsoid. Hence, it is clear that the maximum deviation of the ellipsoid $\Mps \cap \Hplan_\vech$ from its center $\cel$ occurs in the direction $\intbasis_\maxzsv$ (resp. $\kbrace{\Spsbasis_j^*}_{j=\maxzsv+1}^\dimSps$) and is equal to $\svdecr_\maxzsv^{-1}(\Spswidth^2 - \sum_{j=\maxzsv+1}^\dimSo \scap[\Sobasis^*_j]{\vech}^2)$ (resp. $+\infty$) when $\maxzsv=\dimSps$ (resp. $\maxzsv<\dimSps$).



Exploiting \eqref{eq:defMpcapHplan}-\eqref{eq:constrainttot0}, one can randomly sample $\Mps \cap \Hplan_\vech$ by using the procedure described in Algorithm \ref{alg:Sampling1}. The actual distribution according to which the samples are drawn depends on the choice of the distribution on $\pi$ in the first step of the procedure. For example, in the finite dimensional setting, if 
\begin{align}
\pi = \frac{\sum_{j=1}^{\maxzsv-\maxosv} \xi_j^2}{\sum_{j=1}^{\maxzsv-\maxosv+\dimSpspcapSop} \xi_j^2},
\end{align}
where the $\xi_j$'s are independent realizations of a zero-mean Gaussian distribution, the procedure described in Algorithm \ref{alg:Sampling1} draws snapshots from a uniform distribution supported on $\Mps \cap \Hplan_\vech$.  Other choices for the distribution on $\pi$ enable to put more emphasis on some parts of $\Mps \cap \Hplan_\vech$, \ie to draw samples in some regions of the ellipsoid with higher probability (we discuss this option in Section \ref{sec:complexity}). We note however that, as long as the support of the distribution on $\pi$ is $[0,1]$, the samples are drawn from a distribution whose support is equal to $\Mps \cap \Hplan_\vech$.

\begin{algorithm}[t]
\begin{algorithmic}[0]
\State \textbf{inputs}: $\kbrace{\scap[\Sobasis_j]{\vech} }_{j=1}^\dimSo$, $\kbrace{\Spswidth_j}_{j=1}^\nsubspace$, $\kbrace{\Sps_j}_{j=1}^\nsubspace$, $j^*\in\kbrace{1,\ldots,\nsubspace}$
\While{$i< \nsamples$}
	\State 1) Draw $\vech_i$ uniformly at random in the ellipsoid 
	\begin{align}
	\kbrace{\vechtrial : \dist[\vechtrial]{\Sps_{j^*}}\leq \Spswidth_{j^*}}\cap\Hplan_\vech,\nonumber
	\end{align}
	\State \ \ by using the procedure described in the inner loop of Algorithm \ref{alg:Sampling1}. \vspace{0.1cm}
	\State 2)	Acceptance-Rejection:\vspace{0.1cm}
	\If {$\|\projector[\Sps_j^\perp](\vech_i)\| \leq \Spswidth_{j}$ for all $j$}
		\State Include $\vech_i$ to the set of snapshots
	\Else
		\State Reject $\vech_i$
	\EndIf
\EndWhile
\State \textbf{ouput}: $\kbrace{\vech_i \in \Mps \cap \Hplan_\vech}_{i=1}^{\nsamples}$
\end{algorithmic}
\caption{Random Sampling of $\Mps \cap \Hplan_\vech$ for $\nsubspace>1$\label{alg:Sampling2}}
\end{algorithm}

We now focus on the general case where $\Mps$ is defined by \eqref{eq:PriorUoS0} with $\nsubspace>1$. 
In this case $\Mps \cap \Hplan_\vech$ does not usually have any ``desirable'' structure. We note however that the intersection between $\Hplan_\vech$ and each of the $\nsubspace$ degenerate ellipsoids defining $\Mps$ (taken separately) forms an ellipsoid.~We thus propose the ``acceptance-rejection'' strategy described in Algorithm~\ref{alg:Sampling2}.~We consider as a reference ellipsoid, the ellipsoid defined by the intersection of $\Hplan_\vech$ and the $j^*$-th degenerate ellipsoid defining $\Mps$, $j^*\in\kbrace{1,\ldots,\nsubspace}$.~We then draw randomly an element from this ellipsoid by using the procedure described in Algorithm~\ref{alg:Sampling1}. This element, say $\vech_i$, is added to the set of snapshots if it verifies all the constraints defining $\Mps$, \ie 
\begin{align}
\kvvbar{\projector[\Sps_j^\perp](\vech_i)}\leq \Spswidth_{j}\ \mbox{ for all $j\in\kbrace{1,\ldots,\nsubspace}$},\nonumber
\end{align}
and rejected otherwise. 


\subsection{Reduction Based on Point Estimates}\label{sec:point-estimate}


An alternative approach to build a reduced-order model from partial observations may rely on point estimates of the elements of $\Mc$.~More specifically, for each $\vech\in\Ms$, one can compute a point estimate $\hat{\vech}$, obeying some optimality criterion, by combining the partial observations $\kbrace{\scap[\Sobasis_i]{\vech} }_{i=1}^\dimSo$ and some prior information.~Then, an approximation subspace for $\Ms$ can be constructed by considering the manifold of all point estimates, \ie $\hat{\Ms}\triangleq\lbrace\hat{\vech}\rbrace_{\vech\in\Ms}$, as a good surrogate  for $\Ms$. This procedure is summarized in Algorithm~\ref{alg:pointestimate}. 

We note that the procedure described in Algorithm~\ref{alg:pointestimate} shares strong connections with the methodology exposed in \cite{Maday2014PBDW}. In the latter paper, the authors assume that some samples of the solution manifold $\Ms$ are available but some others are only partially observed. They thus propose to complement the set of known samples with point estimates of those only partially observed to compute an approximation subspace for $\Ms$. We see that this setup boils down to the one considered here when \textit{all} the elements of $\Ms$ are partially observed.


~We now make a connection between the procedure described in  Algorithm~\ref{alg:pointestimate} and the material presented in the previous sections.~In particular, we emphasize why such point-estimate procedures may fail in providing reliable results in some situations.

\begin{algorithm}[t]
\begin{algorithmic}[0]
\State \textbf{inputs}: $\kbrace{\kbrace{\scap[\Sobasis_j]{\vech} }_{j=1}^\dimSo}_{\vech\in\Ms}$, $\Mps$ 
\For{$\vech\in\Ms$}
	\State Evaluate $\hat{\vech}$ from $\kbrace{\scap[\Sobasis_j]{\vech} }_{j=1}^\dimSo$ and $\Mps$
\EndFor
\State Apply Algorithm \ref{alg:Greedy} with $\{\hat{\vech}\}_{\vech\in\Ms}$ as input to derive $\Spointa$
\State \textbf{ouput}: $\Spointa$
\end{algorithmic}
\caption{Approximation Subspace from Point Estimates \label{alg:pointestimate}}
\end{algorithm}

There are many choices to compute an estimate $\hat{\vech}$ from partial measurements $\kbrace{\scap[\Sobasis_j]{\vech} }_{j=1}^\dimSo$ and  prior information $\Mps$. Let us consider the following particular option: 
\begin{align}
\hat{\vech} &\in \kargmin_{\vechtrial} \kparen{\sup_{\vech'' \in \Mps \cap \Hplan_\vech} \kvvbar{\vech''-\vechtrial}}. \label{def:wcestimator}
\end{align}
From the worst-case perspective pursued in this paper, \eqref{def:wcestimator} seems indeed to be a sensible choice since $\hat{\vech}$ minimizes the worst-case error over all the elements of $\Sh$ compatible with the received observations and the prior constraints (namely $\Mps \cap \Hplan_\vech$).\footnote{Note that no particular structure for $\Mps$ is assumed in \eqref{def:wcestimator}.} 

 The solution of~\eqref{def:wcestimator} has recently received some attention for some specific choices of $\Mps$ \cite{Binev2015Data}.~When $\Mps$ is defined as the intersection of degenerate ellipsoids as in~\eqref{eq:PriorUoS0}, the authors of \cite{Binev2015Data} showed that $\hat{\vech}$ corresponds to some specific point of $\Mps \cap \Hplan_\vech$ (namely the center of the Chebyshev ball of $\Mps \cap \Hplan_\vech$).\footnote{In fact, as emphasized in \cite[Remark 2.4]{Binev2015Data}, $\hat{\vech}\in\Mps \cap \Hplan_\vech$ as soon as $\Mps \cap \Hplan_\vech$ is a bounded, closed, convex set. }
 In particular, when  $\nsubspace = 1$, $\hat{\vech}$ takes the following simple form
\begin{align}\label{eq:ellipsoidcenter}
\hat{\vech} 
&= 
\sum_{j=1}^\maxzsv \scap[\Sobasis^*_j]{\vech} \svdecr_j^{-1} \Spsbasis^*_j + \sum_{j=\maxzsv+1}^\dimSo \scap[\Sobasis^*_j]{\vech} \Sobasis_j^*,
\end{align}
\ie exactly corresponds to the center of the ellipsoid $\Mps \cap \Hplan_\vech$ defined in \eqref{eq:defMpcapHplan}-\eqref{eq:constrainttot0}.\footnote{When $\maxzsv<\dimSps$, the ellipsoid $\Mps \cap \Hplan_\vech$ is degenerate along the directions $\kbrace{\Spsbasis_{j}}_{j=\maxzsv+1}^\dimSps$. In such a case, there is some ambiguity in the definition of the center of the ellipsoid along these directions. The expression given in \eqref{eq:ellipsoidcenter} then corresponds to the ``center'' with the minimum norm. 
} 

From the perspective of our sampling strategies described in Section \ref{eq:samplingStrategies}, building an approximation subspace from the point estimates $\lbrace\hat{\vech}\rbrace_{\vech\in\Ms}$ is then tantamount to sampling \textit{one} (specific) point of each ellipsoid $\Mps \cap \Hplan_\vech$. This is in constrast with the sampling strategies described in Algorithms \ref{alg:Sampling1} and \ref{alg:Sampling2} where $\nsamples$ points of $\Mps \cap \Hplan_\vech$ are drawn at random.~We may thus expect the point-estimate procedure to lead to performance close to the optimal worst-case solution \eqref{eq:MainProblem} when all the points of $\Mps \cap \Hplan_\vech$ are concentrated around $\hat{\vech}$. When $\nsubspace=1$, this will for example be the case when $\Spswidth$ is small and $\svdecr_i \simeq 1$ for all $i=1,\ldots, \dimSps\leq \dimSo$. On the other hand, when $\svdecr_i \simeq 0$ for some $i$, it is easy to see from \eqref{eq:constrainttot0} that the ellipsoid $\Mps \cap \Hplan_\vech$ will be very elongated along some directions. In such a case, the center of the ellipsoid $\hat{\vech}$ may be a poor representative of the elements of  $\Mps \cap \Hplan_\vech$ and, as a consequence, the approximation subspace computed from $\lbrace\hat{\vech}\rbrace_{\vech\in\Ms}$ may significantly differ from the optimal solution \eqref{eq:MainProblem}. 
 We will illustrate this behavior in our numerical simulations in the next section. We provide below an  illustrative toy example in which the proposed methodology may succeeds in finding a good approximation subspace whereas the methodology based on point estimate \eqref{eq:ellipsoidcenter} is doomed to fail.
 
\begin{example}\label{ex:exampleWC}
We consider the simple case where $\Ms$ is a $\dimSp$-dimensional linear subspace of $\Sh$ and the prior manifold $\Mps$ is defined as in \eqref{eq:PriorSS} for some $\dimSps$-dimensional subspace $\Sps$ and width $\Spswidth$. We assume that by some misfortune 
 the observation subspace $\So$ is orthogonal to $\Ms$. In this case,  we have $\forall \vech\in\Ms$,
\begin{align}
\scap[\Sobasis^*_j]{\vech} = 0  \qquad \forall j=1,\ldots,\dimSo, 
\end{align}
and the ``worst-case optimal'' point estimate defined in \eqref{eq:ellipsoidcenter} is equal to zero. The manifold $\lbrace\hat{\vech}\rbrace_{\vech\in\Ms}$ then reduces to the singleton $\lbrace\0\rbrace$ and a good approximation subspace for $\Ms$ may obviously not be inferred from the latter. 

On the other hand, the posterior manifold $\Mpost$ also takes a simple form when $\So$ is orthogonal to $\Ms$. More specifically, we have\footnote{This can easily be seen by setting $\cel = \0$ in \eqref{eq:defMpcapHplan}-\eqref{eq:constrainttot0} and noticing that $\manifoldperp_\vech$ is equal to $\manifoldperp_0$ $\forall \vech\in\Ms$.}
\begin{align}
\Mpost = \manifoldperp_0,
\end{align}
where $\manifoldperp_0$ corresponds to the ellipsoid defined in \eqref{eq:constrainttot0} with $\vech = \0$. 
 By definition, $\Spost$ will thus be a good approximation subspace for $\Ms$ as soon as $\Ms$ is included in the span of the directions corresponding to the largest variations of $\manifoldperp_0$ (\ie $\kbrace{\Spsbasis_{j}}_{j=\maxzsv+1}^\dimSps$ and the $\intbasis_j$'s associated to the smallest singular values $\svdecr_j$).  We give an illustration of such a scenario (in a slightly more complex setup) in Section \ref{sec:simuresults}. 
 \hfill $\square$
\end{example}

\subsection{Complexity and Implementation Issues}\label{sec:complexity}

In this section, we discuss the computational complexity of Algorithms~\ref{alg:Sampling1} and~\ref{alg:Sampling2}, and elaborate on some related issues. 

Let us first note that  the overall complexity of Algorithms~\ref{alg:Sampling1} and~\ref{alg:Sampling2} is intimately related to the costs of evaluating the sum and the  inner product between two elements of the Hilbert space $\Sh$. In this section, we will assume that these operations can be carried out with a complexity scaling as $\Oc(\dimSh)$. This is for example the case in the finite dimensional setting, which is the one of most interest from an operational point of view. With this convention, we have that the complexity of  Algorithms~\ref{alg:Sampling1} and  \ref{alg:Sampling2} scales as $\Oc\kparen{\dimSo\dimSps\,\dimSh+ \nsamples (\dimSo+\dimSps)\,\dimSh}$ and thus evolves favorably with the problem's dimensions. 

Regarding Algorithm~\ref{alg:Sampling1}, this order of magnitude can be obtained by dividing the analysis of the complexity into the costs of the initial step and the main loop. 
 The most demanding operation in the initial step is the computation of $\kbrace{\Spsbasis_j^*}_{j=1}^\dimSps$, $\kbrace{\Sobasis_j^*}_{j=1}^\dimSo$ and $\kbrace{\svdecr_j}_{j=1}^\mindim$. 
 This task requires to build the Gram matrix associated to the projection operator between the subspaces $\Sps$ and $\So$ and to evaluate its singular value decomposition (see Section \ref{sec:favorable_bases}).
 The first operation involves the computation of $\dimSo\dimSps$ inner products in $\Sh$ and has therefore a complexity scaling as $\Oc\kparen{\dimSo\dimSps\,\dimSh}$. Moreover, the evaluation of the singular value decomposition of an $\dimSo\times\dimSps$ matrix requires at most $\Oc\kparen{\min(\dimSo^2\dimSps, \dimSo\dimSps^2)}$ operations, see \eg \cite[Lecture 31]{Fierro1998Numerical}. Since $\min(\dimSo^2\dimSps, \dimSo\dimSps^2) \leq \dimSo\dimSps\,\dimSh$, the complexity of the initial step is thus of order $\Oc\kparen{\dimSo\dimSps\,\dimSh}$. 

The most demanding operations in the main loop  of Algorithm~\ref{alg:Sampling1} are the evaluations of steps 5 and 6. These tasks can be performed with a complexity scaling at most as $\Oc\kparen{(\dimSo+\dimSps)\,\dimSh}$. Indeed, on the one hand, the uniform sampling of $\vechVpWp$ in $\Sps^\perp\cap\So^\perp$ can be done efficiently as follows: \textit{i)} sample uniformly an element of $\Sh$ (complexity $\Oc\kparen{\dimSh}$); \textit{ii)} set $\vechVpWp = \vech - \projector[\kparen{\Sps^\perp\cap\So^\perp}^\perp](\vech)$ and evaluate $\projector[\kparen{\Sps^\perp\cap\So^\perp}^\perp](\vech)$ with at most $\Oc\kparen{(\dimSo+\dimSps)\,\dimSh}$ operations by noticing that $\kbrace{\Spsbasis_j^*}_{j=1}^\dimSps\cup\kbrace{\intbasis_j}_{j=\maxosv+1}^\maxzsv \cup \kbrace{\Sobasis^*_j}_{j=\maxzsv+1}^\dimSo$
 is an ONB of $\kparen{\Sps^\perp\cap\So^\perp}^\perp$, see Appendix~\ref{sec:charac_MpcapHplan}. On the other hand, the construction of the snapshots $\vech_i$ in step 6 requires the summation of at most $\dimSo+\dimSps+1$ elements of $\Sh$.  Gathering these different elements together, we finally obtain that sampling $\nsamples$ elements of $\Mps \cap \Hplan_\vech$ with Algorithm~\ref{alg:Sampling1} requires a complexity of order $\Oc\kparen{\dimSo\dimSps\,\dimSh+ \nsamples (\dimSo+\dimSps)\,\dimSh}$ as stated above.

The complexity of Algorithm~\ref{alg:Sampling2} is essentially of the same order as the one of Algorithm~\ref{alg:Sampling1} since the latter constitutes the main building block of the former. Nevertheless, we note that the running time of Algorithm~\ref{alg:Sampling2} can be significantly larger than the one of Algorithm~\ref{alg:Sampling1} if the rejection ratio (step 2 in the main loop of Algorithm~\ref{alg:Sampling2}) is important. The choice of the reference ellipsoid $j^*$ appearing in Algorithm \ref{alg:Sampling2} should therefore be made with care in order to decrease as much as possible the rejection ratio. 

Finally, let us note that our study of the complexity is based on the sampling of $\nsamples$ elements of $\Mps \cap \Hplan_\vech$ for \textit{one} $\vech\in\Ms$. In practice, the operations stated in Algorithms~\ref{alg:Sampling1} and \ref{alg:Sampling2} must be repeated for all $\vech\in\Ms$. As a matter of fact, if $\vech\in\Ms$ contains an infinite number of elements, the framework exposed previously cannot, strictly speaking, be applied (it would in particular require to collect a set of observations $\kbrace{\scap[\Sobasis_j]{\vech} }_{j=1}^\dimSo$ for each element in $\Ms$). In such a case, one can nevertheless apply the proposed procedure on a finely-sampled version of the target manifold $\Ms$. All the results discussed previously then carry over by considering the sampled version of $\Ms$ as the new target manifold. One may expect the approximation subspace computed from this sampled manifold to lead to a good approximation subspace for the true manifold $\Ms$ as long as the latter is sampled ``finely enough''. The question of the proper sampling of $\Ms$ is however out of the scope of the present paper and is not further discussed here. We nevertheless mention that, for a given precision, the required number of samples is expected to (typically) scale exponentially with the dimension of the parameter space $\paramSet$. 


To conclude this section, we discuss the choice of the distribution of $\pi$ appearing in step 1 of Algorithm~\ref{alg:Sampling1}. Since any distribution on $\pi$ with support $[0,1]$ defines a distribution on the snapshots supported on $\Mps \cap \Hplan_\vech$, we see that any choice of the distribution on $\pi$ should asymptotically (in the number of samples) leads to the same performance. A proper choice of this  distribution may however have a significant impact on the achievable performance when the number of samples is finite. 
 We advocate below that sampling $\Mps \cap \Hplan_\vech$ with higher probability in the directions of greatest uncertainty may be a good rule of thumbs. 
 
 In order to provide some elements supporting this fact, let us consider the case where $\Mps$ is defined as in \eqref{eq:PriorSS} and $\dimSps=\dimSo$. In such a case, the posterior manifold $\Mpost$ can be simply written as
\begin{align}
\Mpost = \hat{\Ms} \oplus \Ec_0,
\end{align}
where $\hat{\Ms} = \kbrace{\cel}_{\vech\in\Ms}$ with $\cel$ defined in \eqref{eq:cel} and $\Ec_0$ is the ellipsoid defined in \eqref{eq:constrainttot0} with $\vech=\0$.\footnote{$\hat{\Ms}$ is in fact the manifold of the point estimates computed from \eqref{def:wcestimator} with prior \eqref{eq:PriorSS} as discussed in Section \ref{sec:point-estimate}.} 

 Considering this simple expression and assuming that some directions of large variations of the target manifold $\Ms$ are not captured by the point-estimate manifold $\hat{\Ms}$, we see that there is still a hope to identify the latter if they correspond to large variations in $\Ec_0$. In such a case, one may expect the proposed worst-case optimal subspace $\Spost$ to be a better approximation of the target manifold $\Ms$ than $\Spointa$, the subspace evaluated from the sole point-estimate manifold $\hat{\Ms}$. According to this intuition, when only a limited number of snapshots can be drawn from  $\Mps\cap\Hplan$ for each $\vech\in\Ms$, sampling $\Mps\cap\Hplan$ with higher probability in the directions of large uncertainty seems to be a sensible choice.

 
%
%

\section{Simulation Results}\label{sec:simuresults}

\begin{figure}[t!]
\centering
\includegraphics[width=0.55\columnwidth]{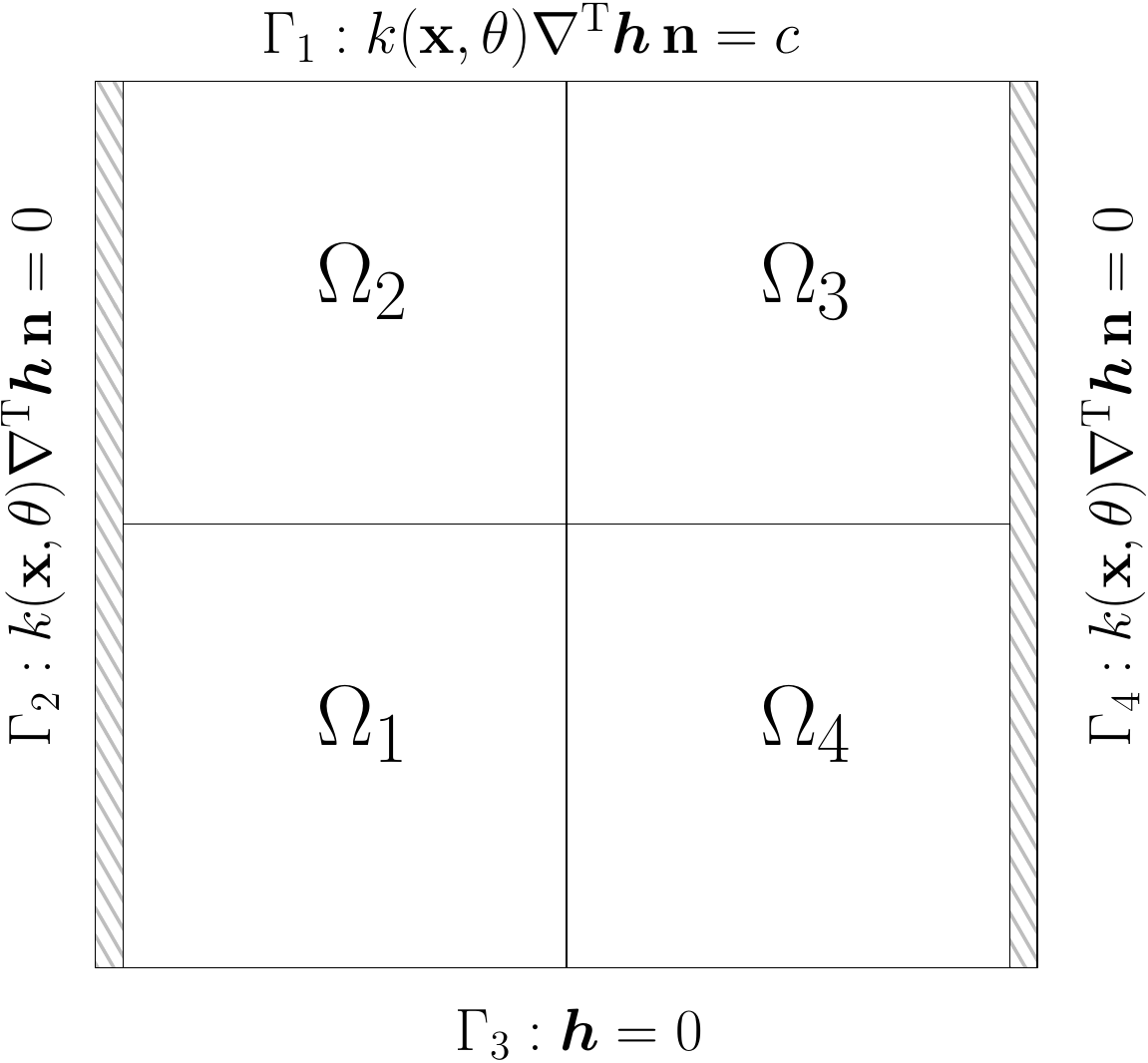}
\caption{\small Schematic representation of the physical system. \label{fig:thermalblockscheme}}
\end{figure}

In this section, we illustrate the performance of the proposed reduction procedures on the standard ``thermal-block'' problem \cite{Patera2006Reduced}: the goal is to evaluate the distribution of the temperature on a plate subject to some boundary conditions, for some specific configurations of the plate's heat conductivity and some external heating source. More specifically, the problem is defined by the following set of differential/boundary equations ($\x\in\Rbb^2$ plays the role of a ``spatial'' variable, $\n$ is a unitary vector normal to the boundary and $\nabla$ is the gradient operator): 
\begin{align}\label{eq:thermalBlockproblem}
\mathrm{PDE}(\vech,\param)=
\left\{
\begin{array}{ll}
\ktranspose{\nabla}\kparen{k(\x,\param)\nabla \vech} = \sourceterm (\x), & \x\in\Omega,\\
k(\x,\param)\ktranspose{\nabla} \vech \,\n =c, & \x\in\Gamma_1,\\
k(\x,\param)\ktranspose{\nabla} \vech \, \n =0, & \x\in\Gamma_2\cup\Gamma_4,\\
\hh =0, & \x\in\Gamma_3,
\end{array}
\right.
\end{align}
where $\Omega = \kbracket{0,1}\times \kbracket{0,1}$ and the boundaries $\Gamma_i$ are defined in Fig.~\ref{fig:thermalblockscheme}.~We assume that the heat conductivity coefficient $k(\x,\theta)$ depends on the parameter $\param = \kbracket{\param_1\, \param_2\,\param_3\,\param_4}\in\paramSet$
 as follows:
\begin{align}
k(\x,\theta) = \sum_{i=1}^4 \param_i\, \mathbb{I}_{\Omega_i}(\x),
\end{align}
where $\mathbb{I}_{\Omega_i}\kparen{\x}$ is the indicator function of $\Omega_i$ and the subdomains $\Omega_i\subset\Omega$ are defined in Fig. \ref{fig:thermalblockscheme}. The definitions of the external heating source $\sourceterm(\x)$ and the boundary parameter $c$ depend on the experiment and are specified below. 

We consider the weak formulation of \eqref{eq:thermalBlockproblem} and approximate its solution via a finite-element method \cite{Brenner2002Mathematical} (see Appendix \ref{app:wformualtion}). The resolution mesh is  chosen fine enough so that the error between the solution of the weak formulation of  \eqref{eq:thermalBlockproblem} and the solution of the finite-element method can be neglected. The discretized system has a dimension equal to $\dimSh=2113$. The solution of the discretized system is computed via the Matlab\textsuperscript{\textregistered} toolbox ``redbKIT'' available at \url{http://redbkit.github.io/redbKIT}.

In our simulations, we apply the subspace identification procedures described in Algorithms \ref{alg:Sampling1}, \ref{alg:Sampling2} and \ref{alg:pointestimate} to the following two setups: 

\paragraph{Setup 1:} We first assume that the target manifold $\Ms$ is defined by the solutions of (the discretized weak formulation of) \eqref{eq:thermalBlockproblem} for 
\begin{align}
\paramSet=\Biggl\lbrace{\param\in\Rbb^4:
\begin{array}{ll}
 \param_i = \param_{\mathrm{min}} + \paramstep\dstep,\ \dstep\in\kbrace{0,\ldots, \nstep}\\
 \param_1=\param_2, \param_3=\param_4
\end{array}\Biggr\rbrace}, \nonumber
\end{align}
where $\param_{\mathrm{min}}=0.1$, $\paramstep=0.1$, $\nstep=20$. In this setup, we assume that $c=1$ and $\sourceterm(\x)=0$ $\forall \x$.~Moreover, we suppose that, for some reasons, the practitioner does not have a perfect knowledge of $\paramSet$ but only knows that $\paramSet\subseteq  \paramSetp$ where
\begin{align}
\paramSetp=\kbrace{\param\in\Rbb^4:
 \param_i = \param_{\mathrm{min}} + \paramstep\dstep,\ \dstep\in\kbrace{0,\ldots, \nstep},  \forall i \in\kbrace{1,\ldots,4}}. \nonumber
\end{align}
We denote by $\Mrs$ the set of solutions of the discretized weak formulation of \eqref{eq:thermalBlockproblem} for $\param\in\paramSetp$. Clearly, we have $\Ms\subseteq\Mrs$ since $\paramSet\subseteq  \paramSetp$.  
Based on this knowledge we derive a prior manifold $\Mps$ of the form \eqref{eq:PriorUoS0} by following the procedure described in Section \ref{sec:favorable_priors}.  More specifically, we apply Algorithm \ref{alg:Greedy} on the elements of $\Mrs$. 
 Letting 
\begin{align}
\Sapproxa_1 \subset \ldots \subset \Sapproxa_\dimSps\nonumber
\end{align}
be the approximation subspaces produced during the first $\dimSps\geq\nsubspace$ iterations of Algorithm \ref{alg:Greedy},  the subspaces $\kbrace{\Sps_i}_{i=1}^\nsubspace$ and scalars $\kbrace{\Spswidth_i}_{i=1}^\nsubspace$ appearing in \eqref{eq:PriorUoS0} are then specified as follows:\footnote{If $\nsubspace =1$ as in \eqref{eq:PriorSS}, we simply set $\Sps  = \Sapproxa_\dimSps$. }
\begin{align} \nonumber
\begin{array}{lcl}
\Sps_1 &= &\Sapproxa_1,\\
            &\vdots& \\
\Sps_{\nsubspace-1} &= &\Sapproxa_{\nsubspace-1},\\
\Sps_\nsubspace &= &\Sapproxa_\dimSps, 
\end{array}
\end{align}
and, $\forall j\in\kbrace{1,\ldots,L}$, 
\begin{align} \nonumber 
\Spswidth_j &= \sup_{\vech\in\Mrs} \dist[\vech]{\Sps_j}. 
\end{align}
We note that \eqref{eq:constraintVeps} is verified by definition of $\Spswidth_j$, 
 so that our working hypothesis \eqref{eq:constraints_prior} is satisfied. In our simulation, we consider the case where $\nsubspace = 1$ and $\nsubspace = 21$. 
 The ONB defining the observation subspace $\So$ is chosen uniformly at random.

 \paragraph{Setup 2:} In this setup, we consider a scenario where the main directions of $\Ms$ are poorly aligned with the observation subspace $\So$, resulting in (potentially) adverse operating conditions for the point-estimate procedure described in Section \ref{sec:point-estimate}. 
 
  In order to precisely describe the experimental setup considered here, we introduce two ONBs, $\kbrace{\Spsbasist_j}_{j=1}^{\dimSpsmax}$ and $\kbrace{\Sobasist_j}_{j=1}^{\dimSpsmax}$, such that 
 \begin{align}
\scap[\Sobasist_i]{\Spsbasist_j} = 0\quad \mbox{ for $i\neq j$},\label{eq:constbasist1}
\end{align}
and 
\begin{align}
\scap[\Sobasist_j]{\Spsbasist_j} = \left\{
\begin{array}{ll}
\corrVW & \mbox{$j= 1 \ldots \dimSp$}, \label{eq:constbasist2}\\
1 & \mbox{otherwise,}
\end{array}
\right.
\end{align}
for some scalar $\corrVW\in(0,1)$  and some integers $\dimSp\leq \dimSpsmax$. We also choose these ONBs so that any $\Spsbasist_j$ and $\Sobasist_j$ verify the homogeneous boundary conditions   of \eqref{eq:thermalBlockproblem}. 
 We then let $\Spbasis_j = (1-\delta^2)^{-\frac{1}{2}} (\Spsbasist_j - \corrVW \Sobasist_j)$ for $j = 1 \ldots \dimSp$. 
 
 We construct $\Ms$ as follows: we set $c=0$, the parameters $\theta$ defining the heat conductivity coefficient $k(\x,\theta)$ are fixed to some constant values and we choose the source term $\sourceterm(\x)$ so that:\footnote{We refer the reader to Appendix \ref{app:wformualtion} for more explanations about the reasons why the manifold \eqref{eq:defMsetup2} can be generated by properly defining the source term $\sourceterm(\x)$.}
\begin{align}\label{eq:defMsetup2}
\Ms &= \manifoldperp_{\mathrm{main}} \oplus \manifoldperp_{\mathrm{perturb}},
\end{align}
where
\begin{align}\nonumber
\begin{array}{ll}
\manifoldperp_{\mathrm{main}} & = \kbrace{\vech = \sum_{j=1}^\dimSp \alpha_j \Spbasis_j : \sum_{j=1}^\dimSp \alpha_j^2 \leq \Spwidth^2_{\mathrm{main}}}, \\
\manifoldperp_{\mathrm{perturb}}&= \kbrace{\vech = \corrVW \sum_{j=1}^{\dimSp}   \beta_j \Sobasist_j + \sum_{j=\dimSp+1}^{\dimSpsmax}  \beta_j \Sobasist_j: \sum_{j=1}^{\dimSpsmax} \gamma_j^{2} \beta_j^2\leq \Spwidth^2_{\mathrm{perturb}}}, 
\end{array}
\end{align}
for some $\Spwidth_{\mathrm{main}}\geq 0$, $\Spwidth_{\mathrm{perturb}}\geq 0$ and $\gamma_j\geq 1$. The solution manifold $\Ms$ is thus the direct sum of two ellipsoids,  $\manifoldperp_{\mathrm{main}}$ and $\manifoldperp_{\mathrm{perturb}}$. We will see below that this particular choice for $\Ms$ 
 together with some specific choice of the observation and prior subspaces lead to a difficult problem for the point-estimate approach. 

In particular, we consider a prior manifold $\Mps$ of the form \eqref{eq:PriorUoS0} by making the following choices:\footnote{If $\nsubspace =1$ as in \eqref{eq:PriorSS}, we simply set $\Sps  = \spa[\kbrace{\Spsbasist_j}_{j=1}^{\dimSps}]$. }
\begin{align} \nonumber
\begin{array}{lcl}
\Sps_1 &= &\spa[\kbrace{\Spsbasist_1}],\\
            &\vdots& \\
\Sps_{\nsubspace-1} &= &\spa[\kbrace{\Spsbasist_j}_{j=1}^{\nsubspace-1}],\\
\Sps_\nsubspace &= &\spa[\kbrace{\Spsbasist_j}_{j=1}^{\dimSps}], 
\end{array}
\end{align}
for some $\dimSps\leq \dimSpsmax$ and, $\forall j\in\kbrace{1,\ldots,L}$, 
\begin{align} \nonumber 
\Spswidth_j &= \sup_{\vech\in\Ms} \dist[\vech]{\Sps_j}. 
\end{align}
The observation subspace is defined as $\So = \spa[\kbrace{\Sobasist_j}_{j=1}^{\dimSo}]$ for some $\dimSo\leq \dimSpsmax$. We note that with this particular choice for $\Sps$ and $\So$, we have from \eqref{eq:constbasist1}-\eqref{eq:constbasist2} that the singular vectors $\{\Spsbasis_j^*\}_{j=1}^\dimSps$, $\{\Sobasis_j^*\}_{j=1}^\dimSo$ simply correspond to $\{\tilde{\Spsbasis}_j\}_{j=1}^\dimSps$, $\{\tilde{\Sobasis}_j^*\}_{j=1}^\dimSo$ whereas the  associated singular values are defined by \eqref{eq:constbasist2}. 

In this particular setup, the main variations of $\Ms$  occur in the direction of $\spa[\kbrace{\Spbasis_j}_{j=1}^\dimSp]$ as soon as  $\Spwidth_{\mathrm{perturb}}\ll \Spwidth_{\mathrm{main}}$. Moreover, if we choose $\corrVW\simeq 0$, we have from \eqref{eq:constbasist1}-\eqref{eq:constbasist2} that $\scap[\Spbasis_i]{\Sobasist_j}\simeq 0$, so that only a small portion of the variations of $\Ms$ in the directions of $\spa[\kbrace{\Spbasis_j}_{j=1}^{\dimSp}]$ can be captured by the observation operator. In particular, if we consider the case where $\corrVW=0$ and $\Spwidth_{\mathrm{perturb}}=0$, we obtain the same setup as in Example \ref{ex:exampleWC}: $\So$ is orthogonal to $\Ms =  \manifoldperp_{\mathrm{main}} \subseteq \spa[\kbrace{\Spbasis_j}_{j=1}^\dimSp]$ and the point-estimate manifold $\hat{\Ms}$ reduces to the singleton $\kbrace{\0}$; on the other hand, we have $\Spbasis_j = \Spsbasis_j^*$ and $\scap[\Sobasis_j^*]{\Spsbasis_j^*}=0$ so that the directions $\{ \Spbasis_j\}_{j=1}^\dimSp$ correspond to directions of large uncertainties in $\Mpost$.\footnote{In particular, as long as prior \eqref{eq:PriorSS} is considered, $\{ \Spbasis_j\}_{j=1}^\dimSp$ corresponds to some degenerate directions 
of the ellipsoid $\manifoldperp_\vech $ in \eqref{eq:constrainttot0}.}

Hereafter, we consider the following scenario: we set $\dimSp = 5$, $\dimSpsmax=50$, $\corrVW = 10^{-4}$, $\Spwidth_{\mathrm{main}}=1$ and $\Spwidth_{\mathrm{perturb}}=10^{-3}$; the weights $\gamma_j$ in the definition of $\manifoldperp_{\mathrm{perturb}}$ are defined as follows:
\begin{align}\label{eq:defgamma}
\gamma_j = \left\{
\begin{array}{ll}
0.85^{-\dimSpsmax} & j = 1\ldots \dimSp,\\
0.85^{-(j-\dimSp)} & j > \dimSp.
\end{array}
\right.
\end{align}
As we will see below, this choice of $\gamma_j$ (together with the definition of $\So$ and $\Sps$) is such that the point-estimate manifold $\hat{\Ms}$ has small variations in the directions of $\spa[\kbrace{\Spbasis_j}_{j=1}^\dimSp]$. 

In our simulation, we chose the  ONB $\kbrace{\Spsbasist_j}_{j=1}^{\dimSpsmax}$ arbitrarily 
and constructed $\kbrace{\Sobasist_j}_{j=1}^{\dimSpsmax}$ so that \eqref{eq:constbasist1}-\eqref{eq:constbasist2} is satisfied. Regarding the definition of $\Mps$, we considered both the cases where $\nsubspace = 1$ and $\nsubspace = 11$. \\[0.2cm]




We now consider the application of the procedures presented  in Sections \ref{eq:samplingStrategies} and \ref{sec:point-estimate} to the two setups described above.~For the construction of $\Mpost$, we draw $\nsamples=5$ elements randomly from $\Mps\cap\Hplan_\vech$ for each $\vech\in\Ms$ by using the methdologies described in Algorithms \ref{alg:Sampling1} and \ref{alg:Sampling2}. In step 1 of Algorithm~\ref{alg:Sampling1}, we define the distribution on $\pi$ as follows:
\begin{align}\label{eq:defdistrPi}
\pi = \left\{
\begin{array}{ll}
\frac{\sum_{j=1}^{\maxzsv-\maxosv} \xi_j^2}{\sum_{j=1}^{\maxzsv-\maxosv+\dimSpspcapSop} \xi_j^2} & \mbox{with probability 0.1},\\
\frac{10^4 \sum_{j=1}^{\maxzsv-\maxosv} \xi_j^2}{10^4 \sum_{j=1}^{\maxzsv-\maxosv} \xi_j^2+\sum_{j=\maxzsv-\maxosv+1}^{\maxzsv-\maxosv+\dimSpspcapSop} \xi_j^2} & \mbox{with probability 0.9},
\end{array}
\right.
\end{align}
where the $\xi_j$'s are independent realizations of a zero-mean Gaussian distribution. The first row in \eqref{eq:defdistrPi} leads to a uniform sampling of $\Mps\cap\Hplan_\vech$ whereas the second one favors the directions of greatest uncertainty. 
In step 4 of Algorithm~\ref{alg:Sampling1}, we approximate the uniform sampling of $\kbrace{d_j}_{\maxzsv+1}^{\dimSps}$ over  $\Rbb^{\dimSps-\maxzsv}$ by a uniform drawing over $\kbracket{-10,10}^{\dimSps-\maxzsv}$. 

\begin{figure*}
\begin{kfig}[1.03\linewidth]{2}
\includegraphics[scale = 2]{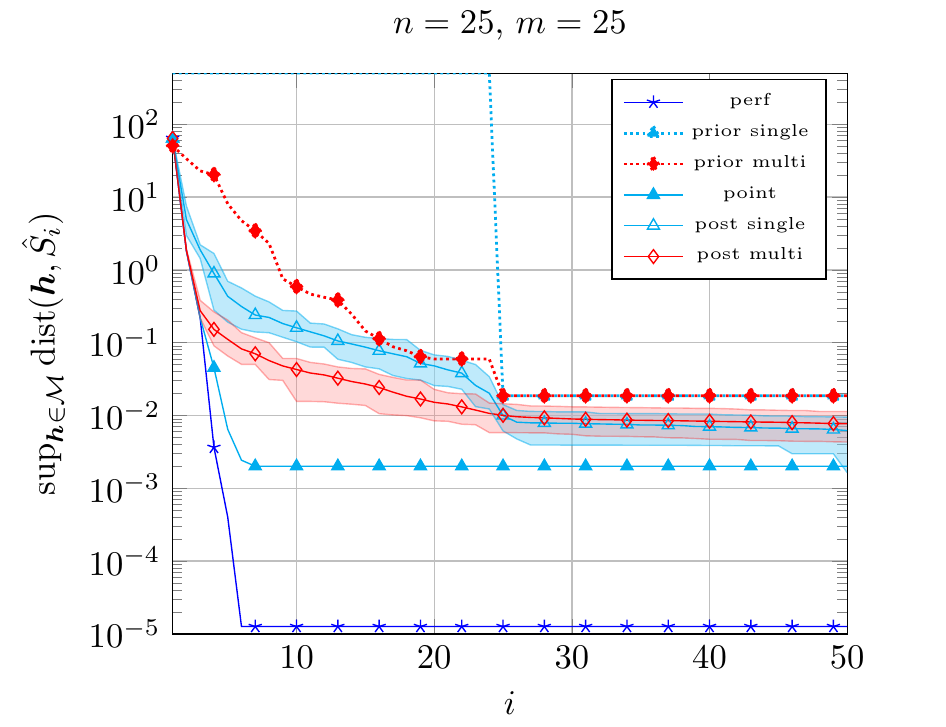}& \includegraphics{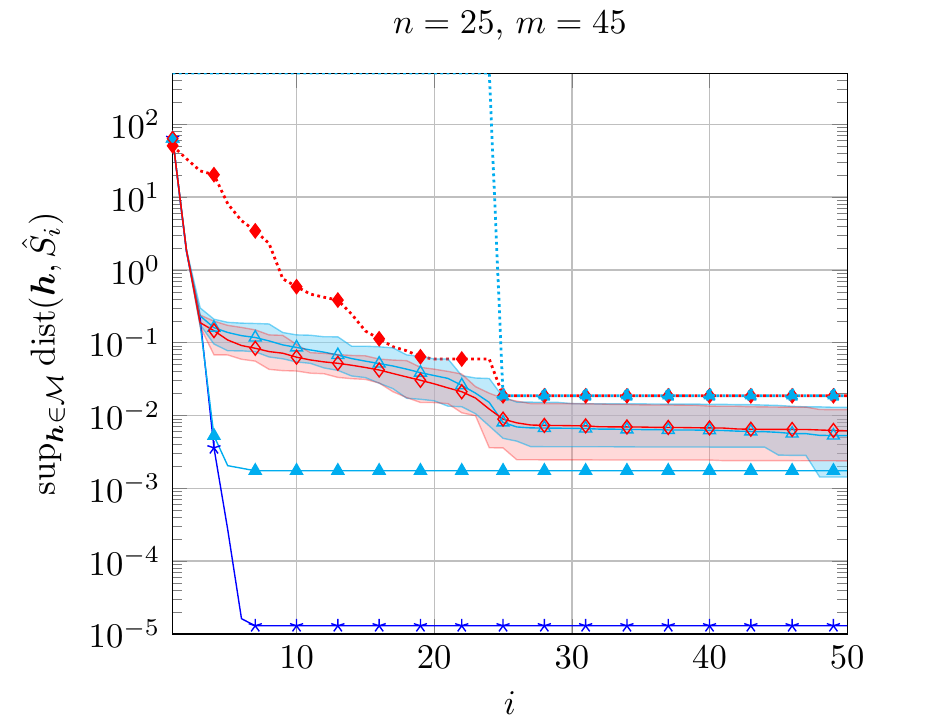}  \kskip
\includegraphics{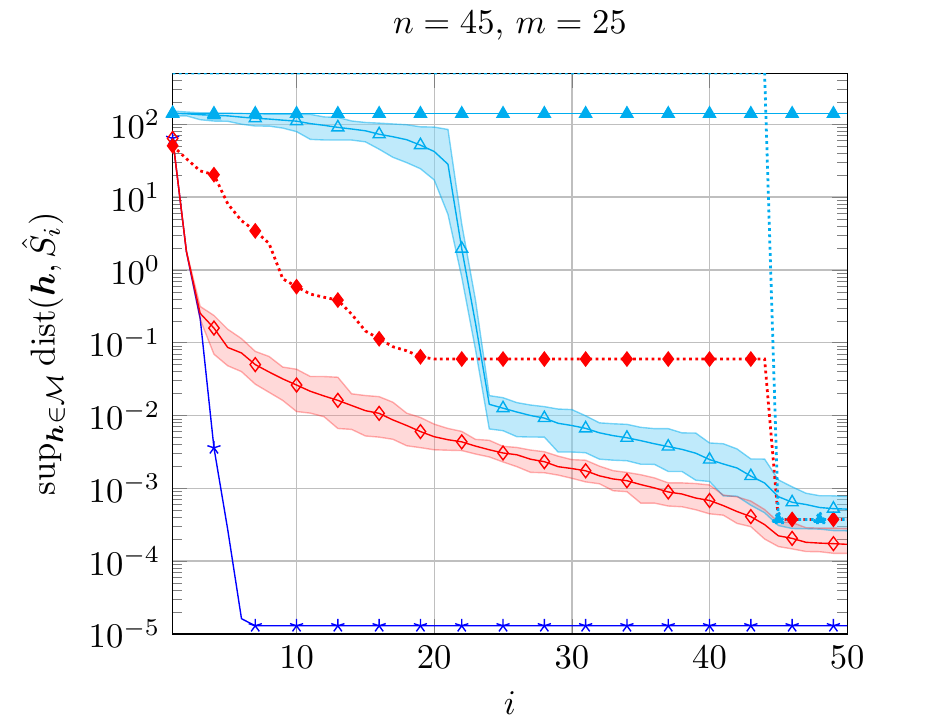}& \includegraphics{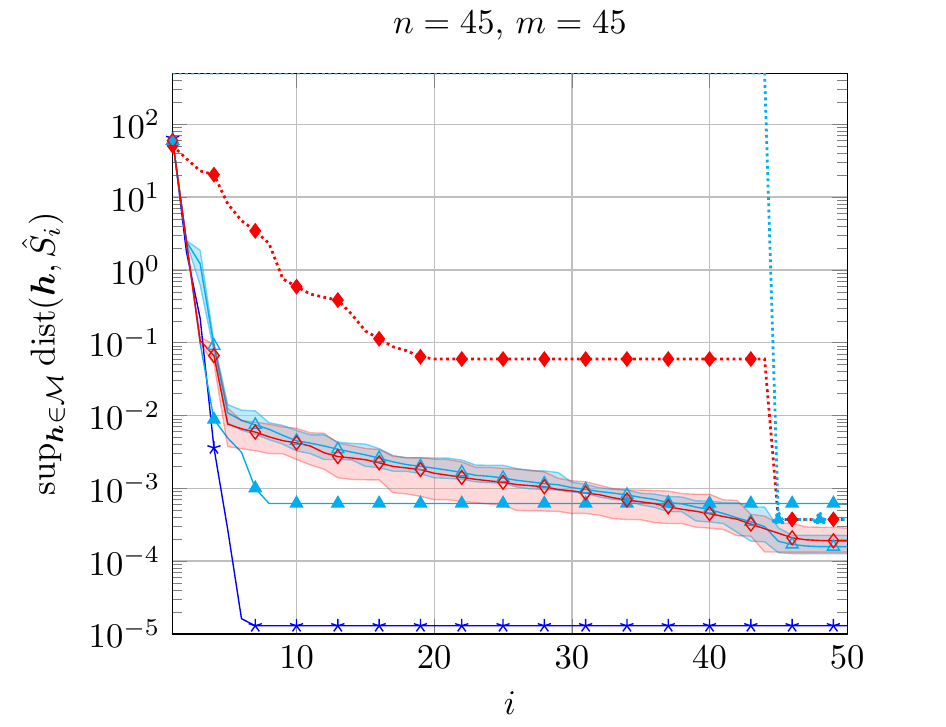} 
\end{kfig}
\caption{\textbf{Setup 1}: Maximum error obtained by projecting  the element of $\Ms$ onto different approximation subspaces. The abscissa represents the dimension $i$ of the approximation subspace. Each figure corresponds to a particular dimension of the prior and observation subspaces, $\Sps$ and $\So$. Regarding the curves ``post single'' and ``post multi'', we ran 20 experiments for each simulation point and averaged the results to obtain the mains curves; the shaded area surrounding these curves represents the interval containing the performance obtained for all the 20 experiments. \label{fig:projerrortruth1}}
\end{figure*}

\begin{figure*}
\begin{kfig}[1.03\linewidth]{2}
\includegraphics{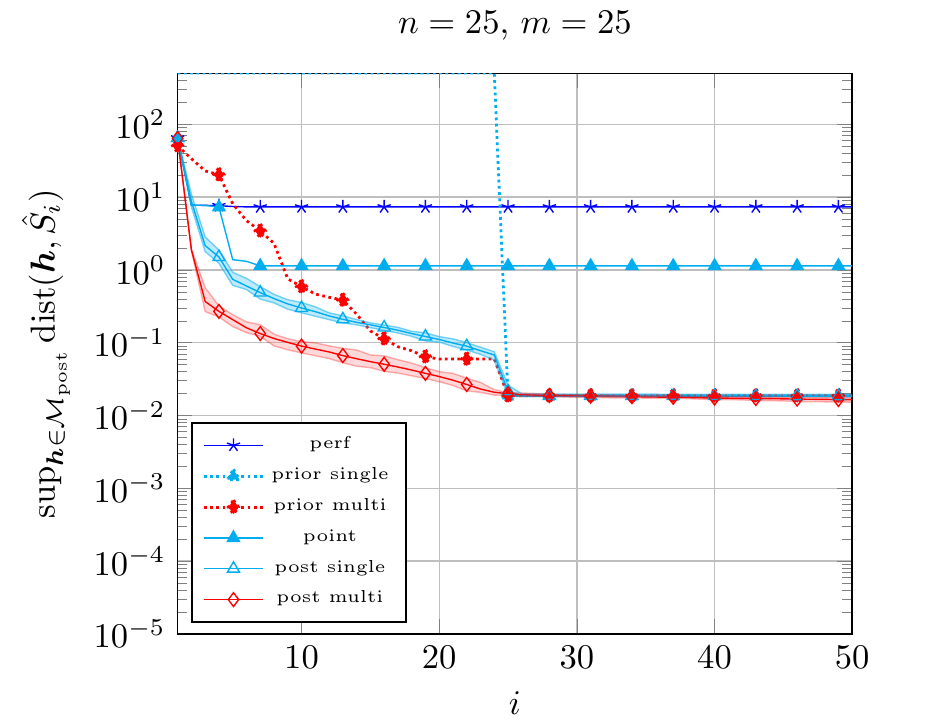}& \includegraphics{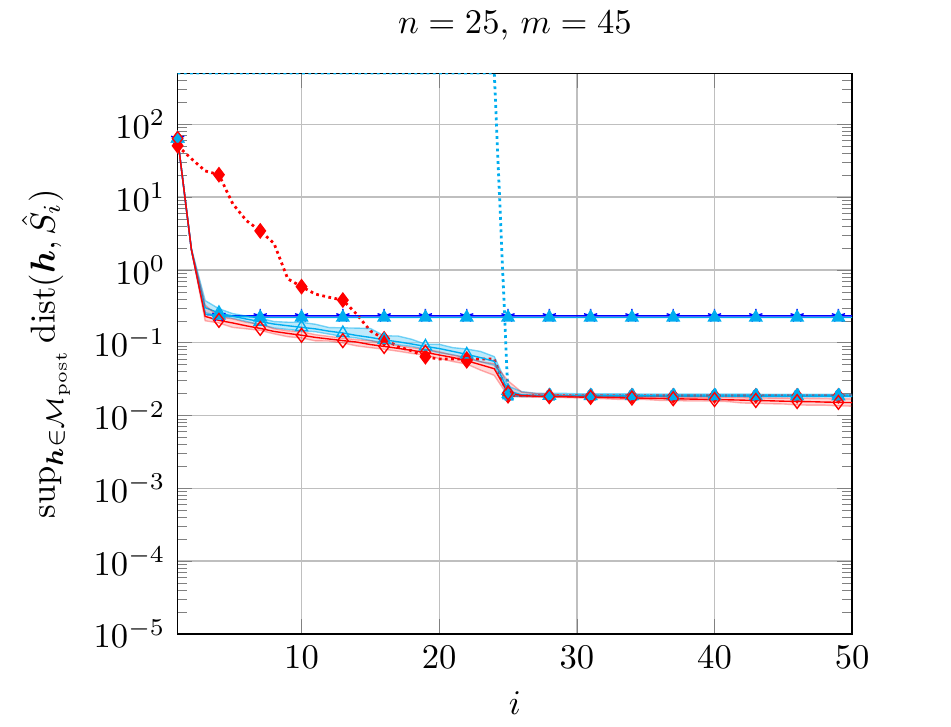}  \kskip
\includegraphics{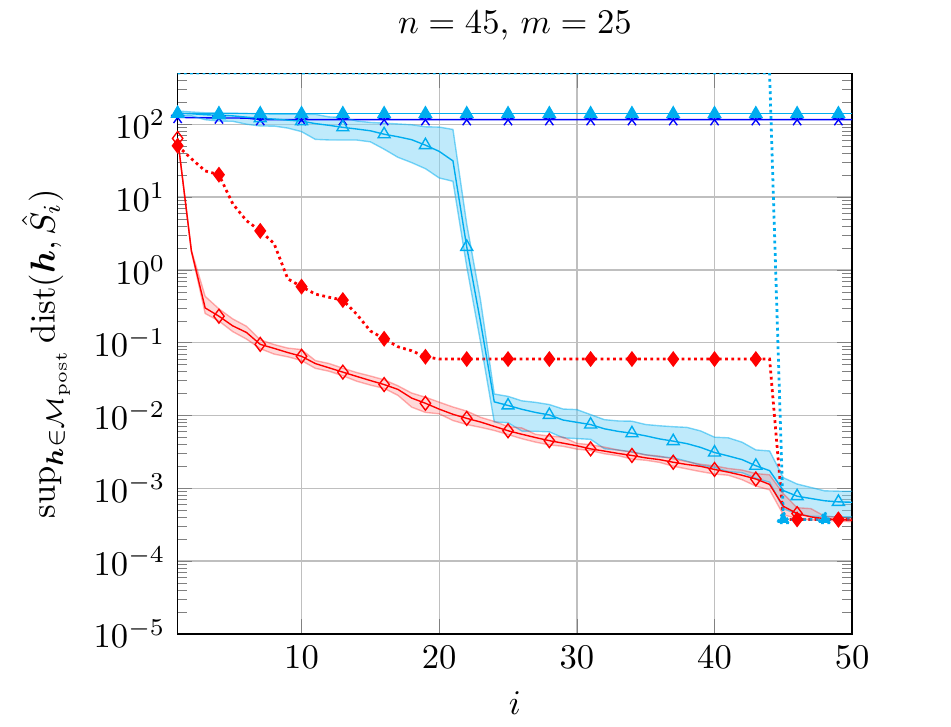}& \includegraphics{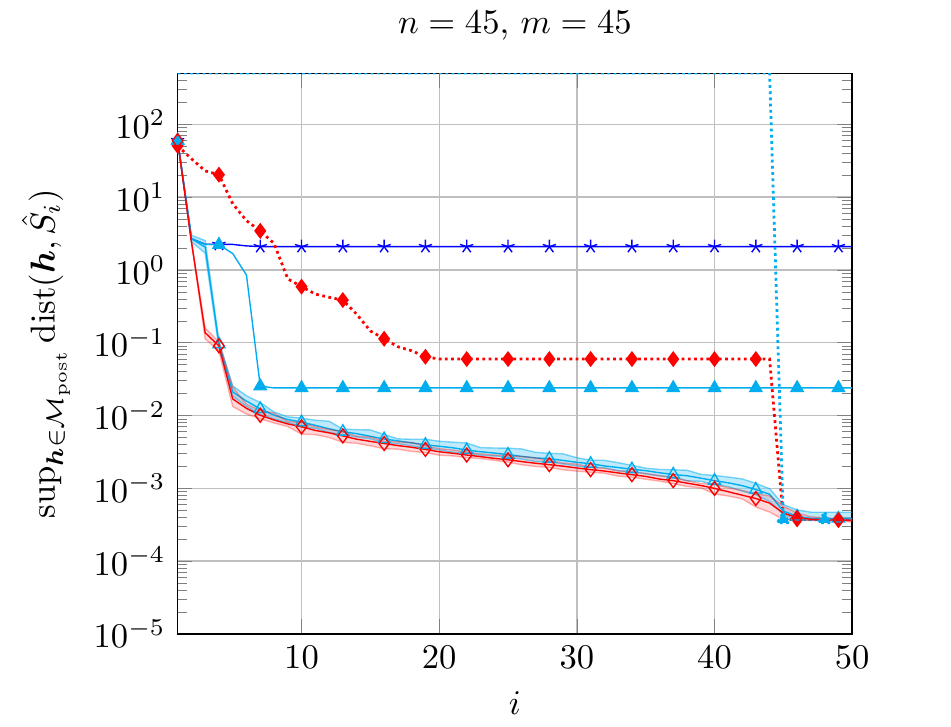} 
\end{kfig}
\caption{\textbf{Setup 1}: Maximum error obtained by projecting  the element of $\Mpost$ onto different approximation subspaces. The abscissa represents the dimension $i$ of the approximation subspace. Each figure corresponds to a particular dimension of the prior and observation subspaces, $\Sps$ and $\So$. Regarding the curves ``post single'' and ``post multi'', we ran 20 experiments for each simulation point and averaged the results to obtain the mains curves; the shaded area surrounding these curves represents the interval containing the performance obtained for all the 20 experiments.\label{fig:projerrorpost1}}
\end{figure*}

\begin{figure*}
\begin{kfig}[1.03\linewidth]{2}
\includegraphics{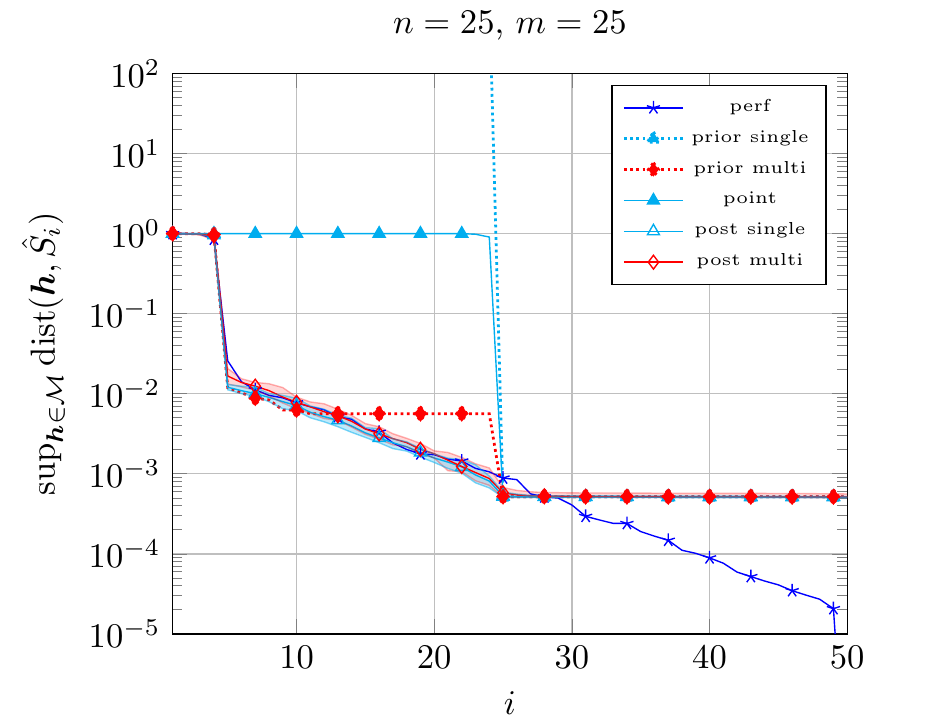}& \includegraphics{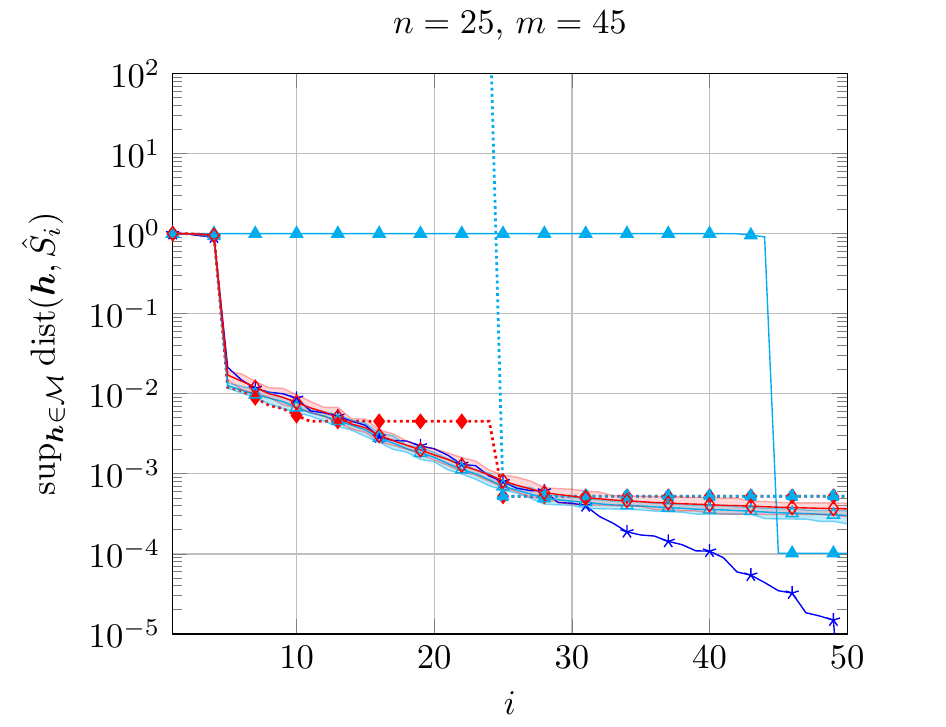}  \kskip
\includegraphics{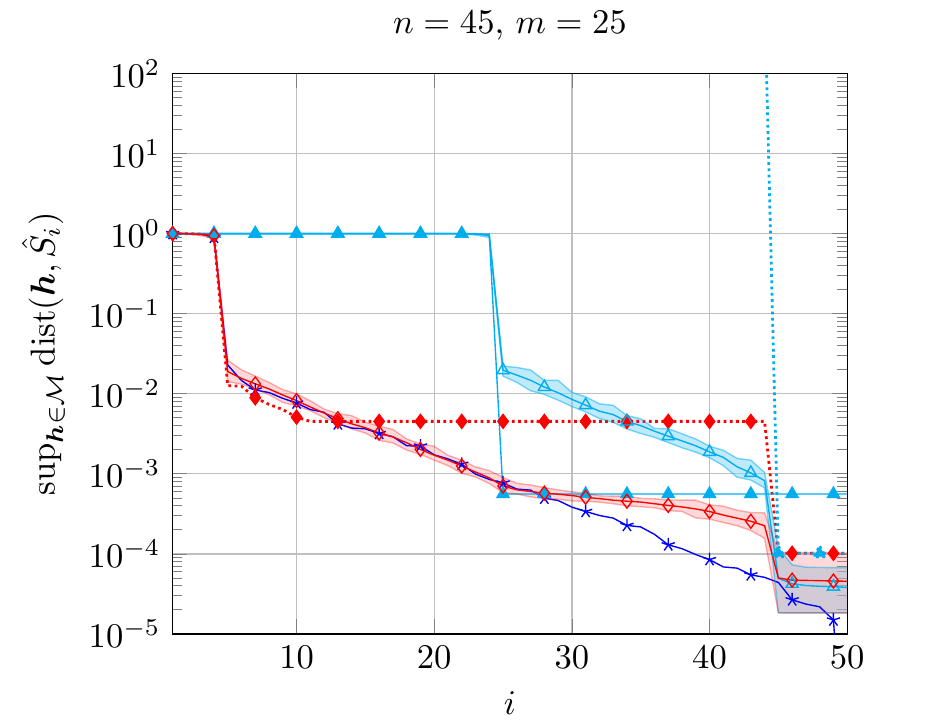}& \includegraphics{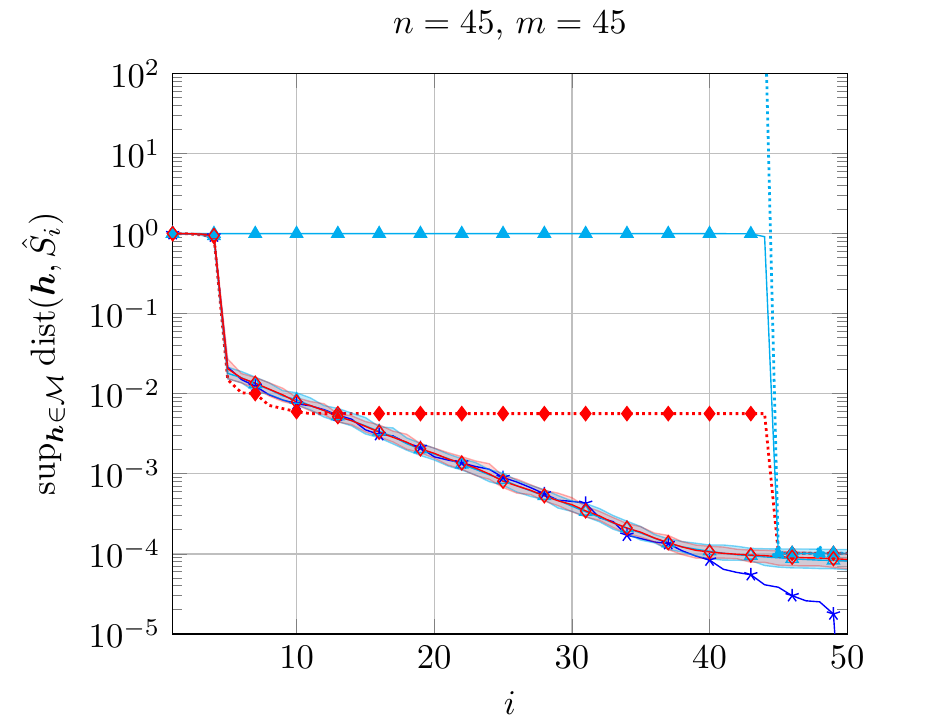} 
\end{kfig}
\caption{\textbf{Setup 2}: Maximum error obtained by projecting  the element of $\Ms$ onto different approximation subspaces. The abscissa represents the dimension $i$ of the approximation subspace. Each figure corresponds to a particular dimension of the prior and observation subspaces, $\Sps$ and $\So$. Regarding the curves ``post single'' and ``post multi'', we ran 20 experiments for each simulation point and averaged the results to obtain the mains curves; the shaded area surrounding these curves represents the interval containing the performance obtained for all the 20 experiments. \label{fig:projerrortruth2}}
\end{figure*}

\begin{figure*}
\begin{kfig}[1.03\linewidth]{2}
\includegraphics{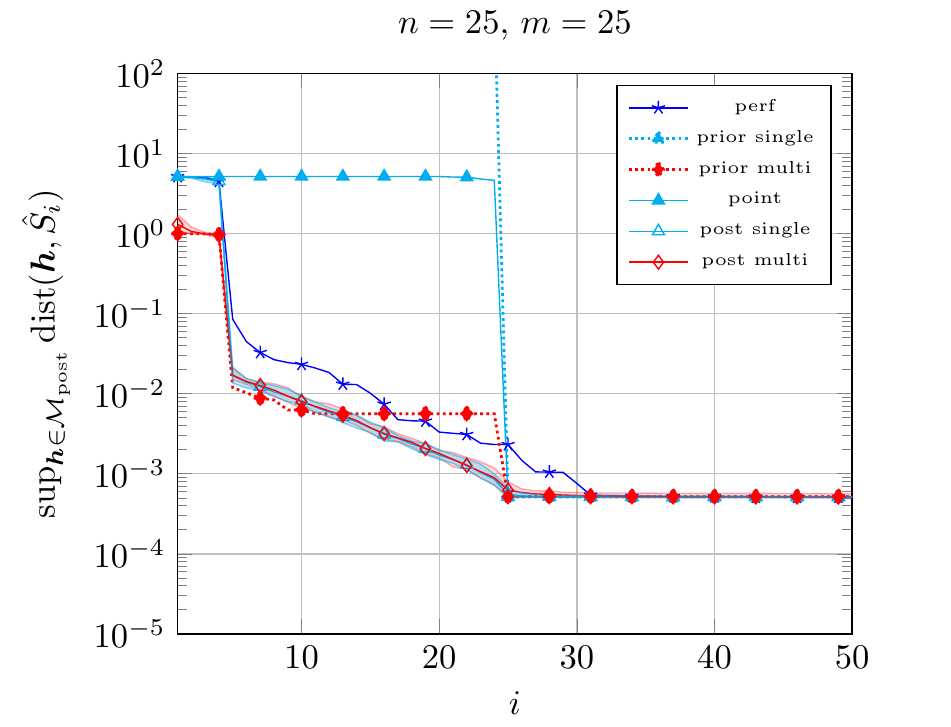}& \includegraphics{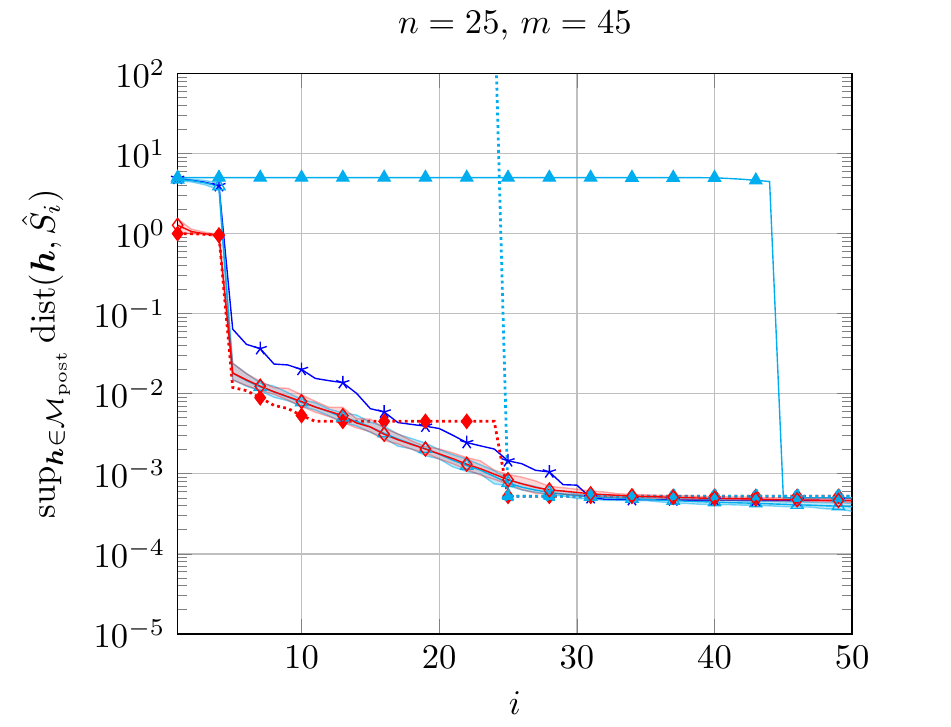}  \kskip
\includegraphics{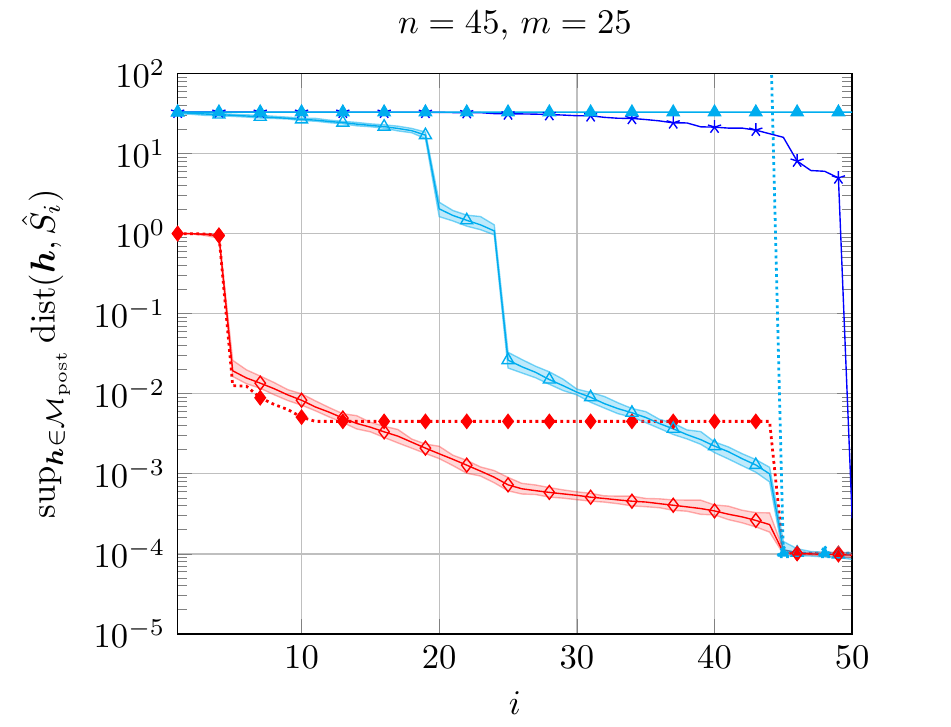}& \includegraphics{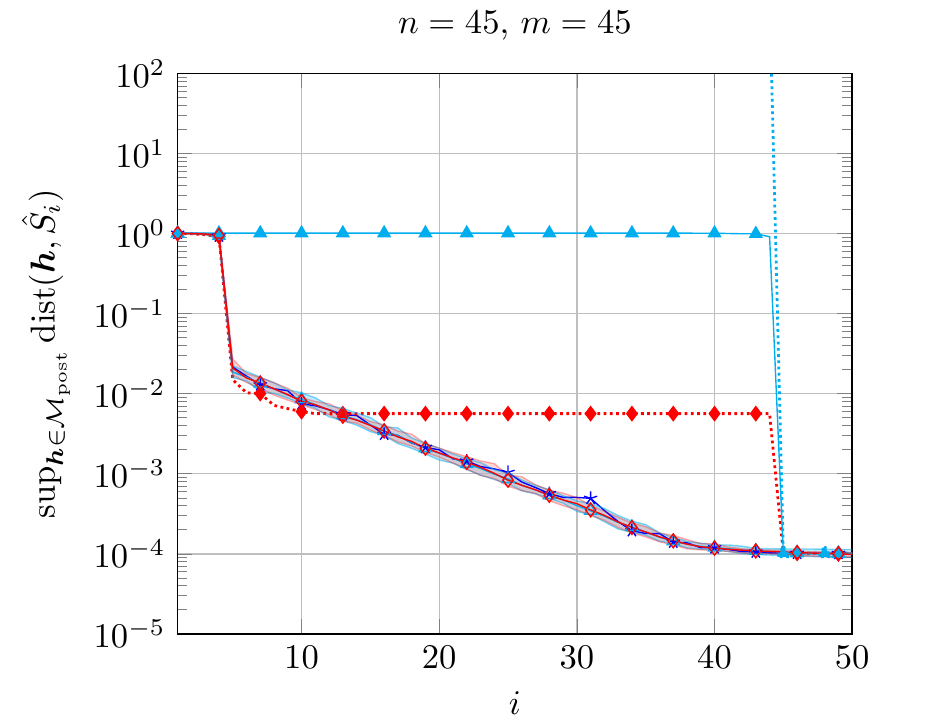} 
\end{kfig}
\caption{\textbf{Setup 2}: Maximum error obtained by projecting  the element of $\Mpost$ onto different approximation subspaces. The abscissa represents the dimension $i$ of the approximation subspace. Each figure corresponds to a particular dimension of the prior and observation subspaces, $\Sps$ and $\So$. Regarding the curves ``post single'' and ``post multi'', we ran 20 experiments for each simulation point and averaged the results to obtain the mains curves; the shaded area surrounding these curves represents the interval containing the performance obtained for all the 20 experiments. \label{fig:projerrorpost2}}
\end{figure*}

Fig.~\ref{fig:projerrortruth1} and~\ref{fig:projerrortruth2} represent the maximum projection error obtained by projecting  the element of $\Ms$ onto different approximation subspaces obtained in Setups 1 and 2, respectively.~These figures thus illustrate the actual reduction performance obtained by different approximation methods. On the other hand, Fig.~\ref{fig:projerrorpost1} and~\ref{fig:projerrorpost2} represent the maximum error obtained by projecting  the elements of $\Mpost$ onto the same approximation subspaces for Setups 1 and 2, respectively. Since $\Mpost$ is the largest manifold compatible with the prior constraint \eqref{eq:constraints_prior} and the received observations (see Lemma \ref{lemma:Mpost}), the curves in Fig.~\ref{fig:projerrorpost1} and~\ref{fig:projerrorpost2} thus provide the worst performance attainable over the set of feasible manifolds.  The parameters $\dimSo$ and $\dimSps$ defining the observation operator and the prior manifold used in each simulation are mentioned on the top of each figure. 

The approximation subspaces considered in all these figures are obtained by applying Algorithm~\ref{alg:Greedy} on:
\begin{itemize}
\item the true manifold $\Ms$ (``perf'').
\item the posterior manifold $\Mpost$ defined from a prior made up of either a single \eqref{eq:PriorSS} or an intersection \eqref{eq:PriorUoS0} of degenerate ellipsoids (respectively ``post single'' and ``post multi'').   
 These curves thus represent the performance achievable by (approximatively) solving the worst-case optimal problem~\eqref{eq:MainProblem}. We note that different experiments may lead to different performance since the snapshots are drawn randomly from $\Mpost$ (see Algorithms~\ref{alg:Sampling1} and~\ref{alg:Sampling2}). This variability  of the results is reported in the figures: for each simulation point, we ran 20 experiments and averaged the results to obtain the main curve; the shaded area surrounding this curve represents the interval containing the performance obtained for all the 20 experiments. 
\item the point-estimate manifold $\hat{\Ms}\triangleq\lbrace\hat{\vech}\rbrace_{\vech\in\Ms}$, where $\hat{\vech}$ is the solution of $\eqref{def:wcestimator}$ with $\Mps$ defined as in \eqref{eq:PriorSS} (``point''). 
 We do not consider the case where $\Mps$ is defined as in \eqref{eq:PriorUoS0} since, in this case,  problem $\eqref{def:wcestimator}$ does not have any simple analytical solution (the problem is actually NP-hard \cite{Binev2015Data}). 
\end{itemize}


We have also reported the Kolmogorov $i$-width, $\Kspec_i(\Mps)$, of the priors used by the different procedures (``prior single'' and ``prior multi'') 
 as a point of comparison. 
%
%
Regarding Fig.~\ref{fig:projerrorpost1} and~\ref{fig:projerrorpost2}, we note that the definition of $\Mpost$ depends on the choice of $\Mps$. The curve ``post multi'' has been computed by using the definition of  $\Mps$ in \eqref{eq:PriorUoS0}. The other curves (``perf'', ``point'' and ``post single'')  have been evaluated by considering the posterior manifold in the case where \eqref{eq:PriorSS} holds. 

We now discuss the performance achieved by the different methodologies mentioned in the paper. The dark blue curve (``perf'') in Fig.~\ref{fig:projerrortruth1} and~\ref{fig:projerrortruth2} corresponds to the performance which can be attained if one has access to the true solution manifold $\Ms$.~The other solid curves illustrate the performance obtained by  exploiting partial observations in the reduction process. 
 As far as the procedures introduced in this paper are concerned (``post single'' and ``post multi''), we observe that the presence of partial observations (almost) always leads to some improvement as compared to the Kolmogorov $i$-width of the prior used in the reduction procedure. This behavior could theoretically be expected from \eqref{eq:upperBKv} since one must have $\sup_{\vech\in\Ms}\dist[\vech]{S_i}\leq \Kspec_i(\Mps)$ as long as the approximation subspace $S_i$ is computed as the solution of \eqref{eq:MainProblem}. Because we only consider an approximate implementation of \eqref{eq:MainProblem} (we draw a finite number of snapshots from  $\Mpost$ and use a greedy strategy to search for a solution of \eqref{eq:MainProblem}), we can nevertheless observe some (limited) degradations of the performance with respect to $\Kspec_i(\Mps)$ in some regions. 
 
 On the other hand, the behavior of the procedure based on point estimates heavily depends on the considered scenario. We see in Fig.~\ref{fig:projerrortruth1} that, as long as Setup 1 is concerned,  the approximation subspace leads to good reduction performance for all choices of $\dimSo$ and $\dimSps$ except for $\dimSo = 25$, $\dimSps = 45$. The latter case corresponds to the scenario where $\Mps$ has more degrees of freedom (\ie $\dimSps = 45$) than the number of observations (\ie $\dimSo = 25$); hence the point estimates \eqref{eq:ellipsoidcenter} computed from the received measurements $\{ \scap[\Sobasis]{\vech}\} $ are poor representatives of the true elements $\vech\in\Ms$. On the other hand, when Setup 2 is considered, we see in Fig.~\ref{fig:projerrortruth2} that the point-estimate approach may, in some cases, not bring any improvement over the performance achievable by the prior manifold ($\dimSo = 25$, $\dimSps = 25$ and $\dimSo = 45$, $\dimSps = 45$) or even degrade the performance ($\dimSo = 45$, $\dimSps = 25$).  
 
 The performance achieved by the point-estimate procedure in Setup 2 can be understood by noticing that the point-estimate manifold  $\hat{\Ms}= \lbrace\hat{\vech}\rbrace_{\vech\in\Ms}$ takes the following form:
\begin{align}
\hat{\Ms} = \kbrace{\hat{\vech} = \sum_{j=1}^{\min(\dimSps,\dimSo)} \beta_j \Spsbasist_j + \sum_{j = \min(\dimSps,\dimSo)+1}^\dimSo \beta_j \Sobasist_j :  \sum_{j=1}^{\dimSo} \gamma_j^2\beta_j^2\leq \Spwidth^2_{\mathrm{perturb}}}.
\end{align}
We first note that $\hat{\Ms}$ is contained in a subspace of dimension $\dimSo$ and the greedy procedure defined in Algorithm~\ref{alg:Greedy} thus necessarily stops after $\dimSo$ iterations. Moreover, since $\corrVW\simeq 0$, the main directions of $\Ms$, namely $\kbrace{\Spbasis_j}_{j=1}^\dimSp$, are such that $\kbrace{\Spbasis_j}_{j=1}^\dimSp \simeq \kbrace{\Spsbasist_j}_{j=1}^\dimSp$. 
 On the other hand, the subspace $\kbrace{\Spsbasist_j}_{j=1}^\dimSp$ corresponds to the directions of smallest amplitudes in $\hat{\Ms}$ because of \eqref{eq:defgamma}. Hence the greedy procedure applied on $\hat{\Ms}$ will select the directions $\kbrace{\Spsbasist_j}_{j=1}^\dimSp$  after only $\dimSo-\dimSp$ iterations, explaining the observed results. In particular, increasing the number of observations may degrade the performance (compare the cases $\dimSo = 25$, $\dimSps = 25$ and $\dimSo = 45$, $\dimSps = 25$). We also note that, unlike in Setup 1,  the point-estimate procedure may  decrease the approximation error in 
the case $\dimSo = 25$, $\dimSps = 45$ for $i>25$. This is due to the fact that, in the particular scenario considered in Setup 2, the $\dimSps-\dimSo=20$ unobserved directions of $\Sps$ are (by construction) orthogonal to the main directions $\kbrace{\Spbasis_j}_{j=1}^\dimSp$ of $\Ms$. 

Regarding the worst-case performance, we see from Fig.~\ref{fig:projerrorpost1} and~\ref{fig:projerrorpost2} that, quite logically, the best results are obtained by the subspaces optimizing the worst projection error over the elements of $\Mpost$. 
In particular, although the point-estimate approach may outperform the proposed procedure regarding the approximation of $\Ms$,  it always leads to inferior performance in the worst-case scenario. 
Interestingly, we also observe that the subspace computed from the true solution manifold $\Ms$ may, in some cases, lead to poor performance, showing that the information contained in the true and the posterior manifolds is quite different.
 We finally note that since $\Ms\subseteq\Mpost$, we have 
\begin{align}
\sup_{\vech\in\Ms} \dist[\vech]{\Sapproxa_i}\leq \sup_{\vech\in\Mpost} \dist[\vech]{\Sapproxa_i},\nonumber
\end{align}
for any approximation subspace $\Sapproxa_i$. Hence, the curves in Fig.~\ref{fig:projerrorpost1} (resp. Fig.~\ref{fig:projerrorpost2})  constitute (approximate\footnote{These bounds are only approximate since, in practice, we  only draw a finite number of snapshots from $\Mpost$.}) upper bounds on those represented in Fig.~\ref{fig:projerrortruth1} (resp. Fig.~\ref{fig:projerrortruth2}). Since these bounds only depend on $\Mpost$, they can always be computed in practice, even if the true solution manifold $\Ms$ is not available.

Regarding the performance achieved by the subspaces computed from the posterior manifold $\Mpost$, the choice of the prior manifold seems to play a crucial role in the minimal projection error which can be attained.~In particular, the smaller the width $\Spswidth_\nsubspace$, the lower the minimal projection error.~The number of observations seems to only have an impact on the minimum subspace dimension required to reach a given approximation error.~For example, in Fig.~\ref{fig:projerrorpost1} the minimal projection error achievable for $\dimSps=25$ is roughly $10^{-2}$ irrespective of the number of observations, whereas an error as low as $\sim 10^{-4}$ can be obtained by setting $\dimSps=45$. In the latter case  we note, as far as the case $\nsubspace=21$ is concerned, that an approximation error of roughly $\sim 10^{-2}$ can be obtained from $i\simeq 10$ when $\dimSo=45$ whereas one needs to increase the dimension of the approximation subspace  up to  $i\simeq 25$ to obtain the same performance when $\dimSo=25$. In the next section, we will provide some theoretical insights into these observations.

\section{Theoretical analysis}\label{sec:thresults}

In this section, we provide a theoretical analysis of the reduction performance achievable within our partially-informed framework. 
 We consider the following simplified scenario:
\begin{align}
\Ms \subseteq  \Mp \cap \Mps,\label{eq:HypMainTh0}
\end{align}
where
\begin{align}\label{eq:HypMainTh}
\begin{array}{rl}
\Mp & \triangleq \{\vech : \dist[\vech]{\Sp}\leq \Spwidth \},\\
\Mps& \triangleq \kbrace{\vech : \dist[\vech]{\Sps}\leq \Spswidth},
\end{array}
\end{align}
for some subspaces $\Sp\subseteq\Sps$, with $\dim(\Sp)=\dimSp$, $\dim(\Sps)=\dimSps$
and some scalars $\Spwidth$, $\Spswidth\geq 0$.~In words, we assume that the unknown manifold $\Ms$ is contained in the intersection of two degenerate ellipsoids  but only one of them (namely $\Mps$) is known a priori.~This scenario thus corresponds to the particular case where $\Mps$ is defined as in \eqref{eq:PriorSS}. 
 
Our goal is to relate the reduction performance obtained by the worst-case optimal approximation subspace $\Spost_i$, defined in \eqref{eq:MainProblem}, to the values of $\dimSp$, $\dimSps$, $\Spwidth$, $\Spswidth$  and the choices of $\Sps$, $\So$. We show below that the latter is closely related to the singular values $\kbrace{\svdecr_j}_{j=1}^\mindim$ of the Gram matrix $\G$ defined in Section \ref{sec:favorable_bases}. 

In order to state our result, we first need to introduce some notations. 
First, let us remind the following notations, introduced in Section \ref{sec:favorable_bases}:  
\begin{align}\nonumber
\begin{array}{ll}
\maxosv & \triangleq \card[ \kbrace{j :\svdecr_j=1}],\\
\maxzsv & \triangleq \card[ \kbrace{j :\svdecr_j>0}].\\
\end{array}
\end{align}
As mentioned previously, the operational meaning of these variables is as follows: $\maxosv$ represents the number of dimensions of $\Sps$ which are included in $\So$, that is $\maxosv=\dim\kparen{\So\cap\Sps}$;  $\dimSps -\maxzsv$ corresponds to the number of dimensions of $\Sps$ which are orthogonal to $\So$, that is $\dimSps-\maxzsv=\dim\kparen{\Sop\cap\Sps}$. In a nutshell, $\maxzsv$ thus represents the number of measurements (out of $\dimSo$) providing information about the position of points in $\Sps$.

We also introduce the following new variable
\begin{align}
\dimSp^*&\triangleq \min\kparen{\dimSps,\dimSp+\dimSps-\maxzsv},\nonumber
\end{align}
and the sequences $\Kspecbound_i$ and $\Kspecboundb_i$:\footnote{We remind the reader that the singular values $\kbrace{\svdecr_j}_{j=1}^\mindim$ are assumed to be sorted in a decreasing order of magnitude.}
\begin{align}\label{eq:UpperBoundKW1}
\Kspecbound_i& \triangleq 
\left\{
\begin{array}{ll}
\infty & i=1,\ldots, \dimSp^*-1,\\
\kparen{\Spwidth+\Spswidth}\, \svdecr_{\maxzsv-(i-\dimSp^*)}^{-1}  & i= \dimSp^*,\ldots, \dimSps-1,\\
\Spswidth & i=\dimSps, \dots, \dimSh-1,
\end{array}
\right.\\
\Kspecboundb_i& \triangleq  \label{eq:UpperBoundKW2}
\left\{
\begin{array}{ll}
\infty & \mbox{$i=1,\ldots, \dimSp+\dim\kparen{\Sop}-1$},\\
\Spwidth & \mbox{$i=\dimSp+\dim\kparen{\Sop}, \dots, \dimSh-1$}.
\end{array}
\right.
\end{align}

 We are now ready to state the following result, whose proof is postponed to  Section \ref{sec:reconstruction_guaran}:
\begin{theorem} \label{th:mainresult}
If $\Ms$ and $\Mps$ verify \eqref{eq:HypMainTh0}-\eqref{eq:HypMainTh}, then the following inequality holds:
\begin{align}\nonumber
\Kspec_i(\Mpost) \leq \min\kparen{\Kspecbound_i,\Kspecboundb_i}\quad \mbox{ for all $i$}. \\[-0.2cm]
\nonumber
\end{align}
\end{theorem}

Interestingly, the upper bounds $\Kspecbound_i$ and $\Kspecboundb_i$ appearing in Theorem~\ref{th:mainresult} only depends on a set of simple parameters defining the partially-informed reduction problem, namely $\dimSp$, $\dimSps$, $\Spwidth$, $\Spswidth$ and the singular values of the Gram matrix $\G$, characterizing the interplay between the prior and observation subspaces, $\Sps$ and $\So$. 
 This few number of parameters enables an easier understanding of the performance achievable by the worst-case reduction methodology presented in  Sections~\ref{sec:WCOMR} and \ref{sec:PracticalImplementation}. 
  Let us recall that the largest projection error induced by the optimal approximation subspace $\Spost_i$ is bounded by the Kolmogorov $i$-width $\Kspec_i(\Mpost)$, see \eqref{eq:upperBKv}.~The upper bound in Theorem \ref{th:mainresult} thus also defines an upper limit on the projection error made by reducing the true, unknown, manifold $\Ms$ in the worst-case optimal subspace $\Spost_i$.

In order to ease the discussion of the result stated in Theorem~\ref{th:mainresult}, we will suppose in the rest of this section that\footnote{We note that satisfying \eqref{eq:HypMainTh0} then requires $\Spwidth\leq\Spswidth$.}
\begin{align}\label{eq:Ms=LDS}
\Mc&= \Mp,
\end{align}
where $\Mp$ is defined in \eqref{eq:HypMainTh}. If $\Spwidth$ is small, \eqref{eq:Ms=LDS} is tantamount to assuming that $\Ms$ has an ``intrinsic dimensionality'' equal to $\dimSp$, \ie the elements of $\Ms$ can vary with no constraint in a $\dimSp$-dimensional subspace $\Sp$ but can only deviate from the latter by a small amount $\Spwidth$ in all the other directions. 

The result stated in Theorem  \ref{th:mainresult} can then be discussed in light of the following comments. 
First, the Kolmogorov $i$-width of $\Ms$ takes a very simple form under hypothesis \eqref{eq:Ms=LDS}, \ie
\begin{align}\label{eq:KwEpsilon} 
\Kspec_i(\Ms)&=
\left\{
\begin{array}{ll}
 \infty & \mbox{if $i<\dimSp$,}\\
 \Spwidth & \mbox{otherwise.}
\end{array}
\right.
\end{align}
The latter provides the best performance which can be achieved by any reduction procedure. In particular, we recall that
\begin{align} \label{eq:lowerprojerror}
\Kspec_i(\Ms) \leq \sup_{\vech\in\Ms} \dist[\vech]{\Spost_i}. 
\end{align}
On the other hand, we have from \eqref{eq:upperBKv} that 
\begin{align}
\sup_{\vech\in\Ms} \dist[\vech]{\Spost_i}\leq \Kspec_i(\Mpost)\leq  \Kspec_i(\Mps). 
\end{align}
As discussed in Section \ref{sec:WCROM},  $\Kspec_i(\Mps)$ is the best-achievable worst-case performance when only $\Mps$ (but no partial observations) is taken into account in the construction of the approximation subspace.~Under assumption \eqref{eq:HypMainTh}, the Kolmogorov $i$-width of $\Mps$ takes again a very simple form, \ie
\begin{align}\label{eq:KwMps}
\Kspec_i(\Mps)&=
\left\{
\begin{array}{ll}
 \infty & \mbox{if $i<\dimSps$,}\\
 \Spswidth & \mbox{otherwise.}
\end{array}
\right.
\end{align}
The gap between $\Kspec_i(\Ms)$ and $\Kspec_i(\Mps)$ represents the potential improvement which can be obtained by the presence of observations. 
  Theorem \ref{th:mainresult} provides an upper bound on the minimal improvement which can be attained by considering $\Spost_i$ as an approximation subspace. In the simple setup considered here, this improvement can be discussed at two different levels:
\begin{itemize}
\item The transition from an infinite to a finite approximation error: we ideally wish to obtain a finite projection error for an approximation subspace whose dimension is as close as possible to $\dimSp$ (we remind the reader that we assume $\Mc= \Mp$ in the present discussion);
\item The approximation error for large dimensions of the approximation subspace: when the size of the approximation subspace increases, we ideally wish to attain a projection error as small as possible.  
\end{itemize}
 Obviously, achieving a finite projection error when reducing the solution manifold $\Ms$ (resp. prior manifold $\Mps$) requires to consider approximation subspaces of dimensionalities greater than or equal to $\dimSp$ (resp. $\dimSps$).~On the other hand, we see from the definitions of $\Kspecbound_i$ and $\Kspecboundb_i$ that the approximation subspace $\Spost_i$ enables to attain a finite projection error as soon as its dimension is greater than or equal to $\min(\dimSp^*,\dimSp+\dim\kparen{\Sop})$. By definition,  $\dimSp^*$ is always greater than or equal to $\dimSp$ but smaller than or equal to $\dimSps$. 
 
The number $\dimSp+\dim\kparen{\Sop}$ corresponds to the ``intrinsic dimensionality'' of $\Ms$ (that is $\dim\kparen{\Sp}\negmedspace~=~\negmedspace\dimSp$) plus the number of dimensions of $\Sh$ which are not measured through our observation operator (that is $\dim\kparen{\Sop}$).~The number $\dimSp^*\triangleq \min(\dimSps,\dimSp+\dimSps-\maxzsv)$ has also an easy interpretation.~The first term in the minimum  corresponds to the number of directions in which we have a priori an infinite uncertainty about the position of $\Ms$, that is $\dim\kparen{\Sps}=\dimSps$.  
 The second term is equal to the intrinsic dimensionality of $\Ms$ (that is $\dim\kparen{\Sp}=\dimSp$) plus the number of components of the prior subspace $\Sps$ which cannot be measured via our observation operator (that is $\dimSps-\maxzsv$). The ``non-observability'' of some directions of $\Sps$ means that if some elements of $\Ms$ have nonzero components in these particular directions, the collected observations provide no information about their magnitudes. 
  When $\dimSp^*=\dimSp+\dimSps-\maxzsv$, the terms ``$\dimSp$'' and ``$\dimSps-\maxzsv$'' thus have different meanings: $\dimSp$ is the number of directions along which the elements of $\Ms$ \textit{do} have a large variation whereas $\dimSps-\maxzsv$ represents the number of directions along which, given the received observations, $\Ms$ \textit{could}  have a large variation.

 We note that $\dimSp+\dim\kparen{\Sop}\leq \dimSp^*$ only if the number of collected observations is large as compared to the dimension of $\Sh$. In particular, if $\dimSh=\infty$, we have $\dim\kparen{\Sop}=\infty$ so that the transition from infinite to finite approximation error always occurs at $\dimSp^*$, \ie $\Kspecbound_i\leq \Kspecboundb_i=\infty$ $\forall i$. Moreover, we have $\dimSp^*< \dimSps$ as soon as $\dimSp<\maxzsv$, that is when the number of observable components of $\Sps$ is larger than the intrinsic dimensionality of $\Ms$. In particular, in the case $\dimSh=\infty$, a finite approximation error occurs at $\dimSp^*=\dimSp$ when $\maxzsv=\dimSps$, that is all the components of $\Sps$ are observed.

Let us now discuss the projection error which can be attained for approximation subspaces of ``sufficient'' dimensionalities. 
 From \eqref{eq:KwEpsilon}-\eqref{eq:lowerprojerror}, we note that the best achievable performance is lower bounded by $\Spwidth$. This projection error can be attained for $i\geq \dimSp+\dim\kparen{\Sop}$ as suggested by Theorem~\ref{th:mainresult} and the behavior of $ \Kspecboundb_i$ in \eqref{eq:UpperBoundKW2}. Nevertheless, as mentioned previously, $\dim\kparen{\Sop}$ can be very large when the number of collected observations is small as compared to the dimension of the ambient space $\Sh$.~In particular, when $\dimSh=\infty$, the transition of $\Kspecboundb_i$ to $\Spwidth$ never occurs. 

In such a case, as suggested by Theorem \ref{th:mainresult}, the behavior of $\Kspec_i(\Mpost)$ is upper bounded by $\Kspecbound_i$. For large values of $i$, we can see that $\Kspecbound_i$ converges to $\Spswidth$. Hence, as far as the number of observations is small with respect to the dimension of $\Sh$, the best projection error which can be achieved  for approximation subspaces of ``moderate'' dimensions seems to be related to the width of the prior manifold $\Mps$. On the other hand, the rate at which the projection error tends to $\Spswidth$ is connected to the conditioning of $\G$.~Indeed, we see from \eqref{eq:UpperBoundKW1} that $\Kspecbound_i$ behaves like $\kparen{\Spwidth+\Spswidth}\, \svdecr_{\maxzsv-(i-\dimSp^*)}^{-1}$ in the range $i= \dimSp^*,\ldots, \dimSps-1$. If $\Spwidth\ll\Spswidth$, we thus have $\Kspecbound_i\simeq \Spswidth\, \svdecr_{\maxzsv-(i-\dimSp^*)}^{-1}\simeq \Spswidth$ provided that $\svdecr_{\maxzsv-(i-\dimSp^*)}\simeq 1$. In such cases, we may thus hope to attain an approximation error close to $\Spswidth$ for $i< \dimSps$, hence improving over the best performance achievable from the prior model, \ie \eqref{eq:KwMps}. 

We remind the reader that the Gram matrix $\G$ (and thus its singular values $\kbrace{\svdecr_j}_{j=1}^\mindim$) characterizes the interplay between the prior and observation subspaces, $\Sps$ and $\So$. In particular, some of the singular values $\svdecr_j$ will be close to one (resp. zero) if some directions of $\Sps$ are almost included in (resp. orthogonal to) $\So$. For example, 
 if $\Sps\subseteq \So$, we obtain $\maxosv=\maxzsv=\dimSps$ and
\begin{align}\nonumber
\begin{array}{ll}
\svdecr_j =1 & j= 1, \ldots, \dimSps, 
\end{array}
\end{align}
so that $\Kspecbound_i= \Spswidth+\Spwidth \simeq \Spswidth$ for $i\geq \dimSp$ when $\Spwidth\ll\Spswidth$. On the other extreme, assuming that $\So$ is orthogonal to $\Sps$ leads to 
\begin{align}\nonumber
\begin{array}{ll}
\svdecr_j =0 & j= 1, \ldots, \mindim, 
\end{array}
\end{align}
so that $\maxosv=\maxzsv=0$ and $\dimSp^*=\dimSps$. In such a case, the performance cannot be improved over that obtained by exploiting the prior manifold $\Mps$ only. \\

To conclude this section, we revisit the empirical results presented in Section \ref{sec:simuresults} in the light of Theorem \ref{th:mainresult}. We focus our discussion on Fig.~\ref{fig:projerrorpost1} which is the more amenable to a connection with our theoretical results. 
%
 First, let us draw a link between the simulation setup in Fig.~\ref{fig:projerrorpost1}  and the parameters defining $\Kspecbound_i$ and $\Kspecboundb_i$.  
 From the dark blue plain curve in Fig.~\ref{fig:projerrortruth1}, it can be seen that $\Ms$ is included in a degenerate ellipsoid with $\dimSp=4$ and width $\Spwidth\sim10^{-5}$.\footnote{In the discussion below, we implicitly assume that the 4-dimensional subspace defining the ellipsoid is included in the prior subspace $\Sps$, so that our simulation results can be discussed in light of Theorem \ref{th:mainresult}.}~As far as the case $\nsubspace=1$ was concerned, we considered two different dimensions for the subspace $\Sps$, that is $\dimSps=25$ or $\dimSps=45$.~These two choices led to $\Spswidth\sim10^{-2}$ and $\Spswidth\sim10^{-4}$, respectively. Hence, we have $\Spwidth\ll\Spswidth$ in both cases. The observation subspace $\So$ was drawn uniformly at random, so that the situation $\maxosv=0$ and $\maxzsv=\min(\dimSo,\dimSps)$  occurs almost surely in our simulations; in this case,  $\dimSp^*$ takes the simple form $\dimSp^*=\dimSp+\dimSps-\mindim$. We considered a number of observations equal to either $\dimSo=25$ or $\dimSo=45$. 
 
 For these choices of parameters and for the dimensions of the approximation subspaces considered in our simulation (that is $i=1,\ldots, 50$), the bound $\Kspecbound_i$ prevails in Theorem \ref{th:mainresult}. Indeed, $\dim(\Sop)\geq \dimSh-45=2113-45=2068$ so that $\dimSp+\dim(\Sop)\geq 2072$ and $\Kspecboundb_i=\infty$ for the range of interest. 
 
 Although Theorem \ref{th:mainresult} only provides an upper bound on $\Kspec_i(\Mpost)$,  the empirical performance presented in Fig.~\ref{fig:projerrorpost1} seems to be in good accordance with the latter.  
 First, we can notice that the performance achieved by the worst-case approximation subspace $\Spost_i$ (plain magenta curve) saturates at $\Spswidth$ for large $i$. As predicted by Theorem \ref{th:mainresult}, the number of observations does not have an effect on this error floor but rather impacts the rate at which the latter is reached. 
 Theorem \ref{th:mainresult} also gives good insights into the range of dimensions in which small approximation errors can be attained. Of particular interest for us is the case where $\dimSo=25$, $\dimSps=45$ for which $\dimSp^*=\dimSp+\dimSps-\dimSo=\dimSp+20$ ($\dimSp^*=\dimSp$ in all other cases); $\dimSps-\dimSo=20$ corresponds to the number of dimensions of $\Sps$ which are not observed. In practice, we can therefore not expect to have a small (worst-case) projection error if $i\leq 20$. This is what we observe\footnote{The error is not infinite in this range  because of the way we generate the posterior manifold $\Mpost$, see Section \ref{sec:simuresults}. } in Fig.~\ref{fig:projerrorpost1}.

\section{Proof of Theorem \ref{th:mainresult}} \label{sec:reconstruction_guaran}

In this section, we provide the main steps of the proof of Theorem \ref{th:mainresult} (the technical details are postponed to the appendices). 
The simple observation underlying the proof is as follows: for any $\Mposts$ such that $\Mpost\subseteq \Mposts$ and  $\forall \Sapprox_i$ with $\dim(\Sapprox_i)=i$, we have
\begin{align}
\Kspec_i(\Mpost)\leq \max_{\vech\in \Mposts} \dist[\vech]{\Sapprox_i}.  \label{eq:relaxation_Kw}
\end{align}
The result stated in Theorem \ref{th:mainresult} then follows from \eqref{eq:relaxation_Kw} and some specific choices for $\Mposts$ and $\Sapprox_i$. These choices are described in Sections \ref{sec:charac_MPcapHplan} and \ref{sec:derivationUpperBound}. 

\subsection{Definition of $\Mposts$ \label{sec:charac_MPcapHplan}}

We give hereafter three possible choices for $\Mposts$ that will be exploited in our proof in Section \ref{sec:derivationUpperBound}.
 First, we have by definition 
\begin{align}
\Mpost\subseteq\Mps, \label{eq:defMposts1}
\end{align}
so that $\Mposts=\Mps$ is a valid choice. 
Moreover, since we assume that $\Ms\subseteq\Mp$, we also have
\begin{align}
\Mpost &\subseteq  \cup_{\vech\in\Mp}\Hplan_\vech,  \label{eq:UBMpsintstep2}\\
\Mpost &\subseteq  \Mps\cap\cup_{\vech\in\Mp}\Hplan_\vech.  \label{eq:UBMpsintstep}
\end{align}
The right-hand sides of \eqref{eq:UBMpsintstep2}-\eqref{eq:UBMpsintstep} thus constitute two other possible choices for $\Mposts$. 

 A precise mathematical characterization of $\cup_{\vech\in\Mp} \kparen{\Mps\cap\Hplan_\vech}$ is given by (see Appendix \ref{app:caractMpost}) 
\begin{align}
\Mps\cap\cup_{\vech\in\Mp}\Hplan_\vech&=\Spext \oplus \manifoldperp, \label{eq:formunionEllipsoid}
\end{align}
where 
\begin{itemize}
\item $\Spext \triangleq \Sp \oplus \spa[ \kbrace{\Spsbasis_j^*}_{j=\maxzsv+1}^\dimSps]$  is a linear subspace 
\item $\manifoldperp$ is an ellipsoid defined as follows: $\vech\in \manifoldperp$ if and only if $\vech$ can be written as 
\begin{align} \label{eq:defmanifoldperp1}
\vech &= \sum_{j=\maxosv+1}^{\maxzsv} \svdecr_j^{-1}\kparen{ \kparen{1-\svdecr_j^2}^{\frac{1}{2}}a_j - b_j} \intbasis_j 
+ \sum_{j= 1}^{\dimSpspcapSop} c_j \orthSoSpsbasis_j + \sum_{j=1}^\dimSo a_j \Sobasis_j^*,
\end{align}
where 
 $\intbasis_j\triangleq \kparen{1-\svdecr_j^2}^{-\frac{1}{2}}\kparen{\Spsbasis^*_j - \svdecr_j \Sobasis^*_j}$, 
$\kbrace{\orthSoSpsbasis_j}_{j=1}^{\dimSpspcapSop}$ is an (arbitrary) ONB of $\Sop \cap \Spsp$, and the coefficients $a_j$, $b_j$, $c_j$  obey the following constraints:
\begin{align}\label{eq:defmanifoldperp2}
\left\{ 
\begin{array}{l}
 \displaystyle{\sum_{j=\maxzsv+1}^\dimSo a_j^2+\sum_{j= \maxosv+1}^\maxzsv b_j^2 +\sum_{j= 1}^{\dimSpspcapSop} c_j^2 \leq \Spswidth^2},\\[0.2cm]
 \displaystyle{\sum_{j=1}^\dimSo a_j^2\leq \Spwidth^2}.
\end{array}
\right.
\end{align}
\end{itemize}
The vectors $\Spsbasis^*_j$, $\Sobasis^*_j$ and $\intbasis_j$ which appear in the above characterization are those introduced in Section \ref{sec:favorable_bases} with the following conventions: $\Sps$ is the subspace characterizing $\Mps$ in \eqref{eq:HypMainTh}; $\So$ is the observation subspace introduced in \eqref{eq:obseq}. 


\subsection{Computation of upper bounds on $\Kspec_i(\Mpost)$} \label{sec:derivationUpperBound}

In this section, we prove the upper bound stated in Theorem~\ref{th:mainresult} by exploiting \eqref{eq:relaxation_Kw} for different choices of $\Mposts$ and $S_i$. The upper bound $\Kspecboundb_i$ straightforwardly derives from the definition of $\Mposts$ given in \eqref{eq:UBMpsintstep2}; $\Kspecbound_i$ results from a combination of the upper bounds obtained from the definition of $\Mposts$ in \eqref{eq:defMposts1} and \eqref{eq:UBMpsintstep}. The detailed calculations are provided below. 

\subsubsection{Case $\Mposts= \cup_{\vech\in\Mp}\Hplan_\vech$\label{sec:defKbarbar}}:
We first note that 
\begin{align}
\cup_{\vech\in\Mp}\Hplan_\vech = \kbrace{\vech : \dist[\vech]{\cup_{\vech\in\Sp}\Hplan_\vech}\leq \Spwidth},
\end{align}
where $\cup_{\vech\in\Sp}\Hplan_\vech$ is a subspace of dimension at most equal to $\dimSp + \dim\kparen{\So^\perp}$. Hence,  if we let $\Sapprox_i$ be such that $\cup_{\vech\in\Sp}\Hplan_\vech \subseteq \Sapprox_i$ for $i\geq \dimSp + \dim\kparen{\So^\perp}$, we obtain
\begin{align}
\Kspec_i(\Mpost) \leq \Spwidth.  
\end{align}

\subsubsection{Case $\Mposts=\Mps$ \label{sec:defKbar1}}: 
Since $\Kspec_i(\Mpost)\leq \Kspec_i(\Mps)$, setting $\Mposts=\Mps$ in \eqref{eq:relaxation_Kw} leads to 
\begin{align}\label{eq:bound1}
\Kspec_i(\Mpost) \leq 
\left\{
\begin{array}{ll}
\infty & \mbox{$i=1,\ldots, \dimSps-1$,}\\
\Spswidth & \mbox{$i=\dimSps,\ldots, \dimSh-1$}.
\end{array}
\right.
\end{align}


\subsubsection{Case $\Mposts=\Mps\cap\cup_{\vech\in\Mp}\Hplan_\vech$\label{sec:defKbar3}}: 
We note that the dimension of the subspace $\Spext$ defined in Section \ref{sec:charac_MPcapHplan} is at most equal to $\min\kparen{\dimSps,\dimSp+\dimSps-\maxzsv}=\dimSp^*$, where   $\dimSp^* \le \dimSps$ follows from the  inclusion $\Sp\subseteq\Sps$. 
For $i =\dimSp^*,\ldots,\dimSp^*+\maxzsv-\maxosv-1$, we set  $\Sapprox_i = \Spext \oplus R$, where $R$ is an $(i-\dimSp^*)$-dimensional subspace given by 
\begin{align}
R \triangleq \spa[\kbrace{\intbasis_j}_{j=\maxzsv-(i-\dimSp^*)+1}^{\maxzsv}].\nonumber
\end{align}
We note that 
\begin{align}
\max_{\vech\in \Mposts} \dist[\vech]{\Sapprox_i} = \max_{\vech\in \Mposts} \kvvbar{\projector[\Sapprox_i^\perp](\vech)}. \nonumber
\end{align}
Moreover, because of the  characterization of $\Mps\cap\cup_{\vech\in\Mp}\Hplan_\vech$ in \eqref{eq:formunionEllipsoid}-\eqref{eq:defmanifoldperp2}, we have that 
\begin{align}
\kvvbar{\projector[\Sapprox_i^\perp](\vech)} 
&= \kvvbar{\projector[\Sapprox_i^\perp](
\sum_{j=\maxosv+1}^{\maxzsv-(i-\dimSp^*)} \svdecr_j^{-1}\kparen{ \kparen{1-\svdecr_j^2}^{\frac{1}{2}}a_j - b_j} \intbasis_j+ \sum_{j= 1}^{\dimSpspcapSop} c_j \orthSoSpsbasis_j + \sum_{j=1}^\dimSo a_j \Sobasis_j^*)},\nonumber\\
&\leq \kvvbar{\sum_{j=\maxosv+1}^{\maxzsv-(i-\dimSp^*)} \svdecr_j^{-1}\kparen{ \kparen{1-\svdecr_j^2}^{\frac{1}{2}}a_j - b_j} \intbasis_j+ \sum_{j= 1}^{\dimSpspcapSop} c_j \orthSoSpsbasis_j + \sum_{j=1}^\dimSo a_j \Sobasis_j^*},\nonumber\\
&= \kparen{\sum_{j=\maxosv+1}^{\maxzsv-(i-\dimSp^*)} \svdecr_j^{-2}\kparen{ \kparen{1-\svdecr_j^2}^{\frac{1}{2}}a_j - b_j}^2+\sum_{j= 1}^{\dimSpspcapSop} c_j^2 + \sum_{j=1}^\dimSo a_j^2}^{\frac{1}{2}},\label{eq:longequation}
\end{align}
holds for any $\vech\in \Mps\cap\cup_{\vech\in\Mp}\Hplan_\vech$, where the 
 last equality follows from the orthogonality of $\kbrace{\Sobasis_j^*}_{j=1}^\dimSo$, $\kbrace{\intbasis_j}_{j=\maxosv+1}^\maxzsv$ and $\kbrace{\orthSoSpsbasis_j}_{j=1}^{\dimSpspcapSop}$ (see Appendix \ref{sec:charac_MpcapHplan}).  

Finally, since the coefficients $a_i$, $b_i$ and $c_i$ must satisfy the constraints in \eqref{eq:defmanifoldperp2}, we have
\begin{align}
\max_{\vech\in \Mposts} \dist[\vech]{\Sapprox_i} 
&\leq \kparen{\svdecr_{\maxzsv-(i-\dimSp^*)}^{-2}\kparen{ \kparen{1-\svdecr_{\maxzsv-(i-\dimSp^*)}^2}^{\frac{1}{2}}\Spwidth + \Spswidth}^2 + \Spwidth^2}^{\frac{1}{2}},\nonumber\\ 
&\leq \svdecr_{\maxzsv-(i-\dimSp^*)}^{-1} \kparen{\Spswidth+\Spwidth}.\label{eq:bound2}
\end{align}
%

 We follow the same reasoning for $i =\dimSp^*+\maxzsv-\maxosv, \ldots, \dimSp^*+\maxzsv-1$. 
  We set
\begin{align}
\Sapprox_i = \Spext \oplus \spa[\kbrace{\intbasis_j}_{j=\maxosv+1}^{\maxzsv}] \oplus R,\nonumber
\end{align}
where $R$ is some arbitrary $(i-(\dimSp^*+\maxzsv-\maxosv))$-dimensional subspace. 
 From the particularization of $\Mps\cap\cup_{\vech\in\Mp}\Hplan_\vech$ in \eqref{eq:formunionEllipsoid}-\eqref{eq:defmanifoldperp2}, we have for any $\vech\in \Mps\cap\cup_{\vech\in\Mp}\Hplan_\vech$:
\begin{align}
\kvvbar{\projector[\Sapprox_i^\perp](\vech)} 
&= \kvvbar{\projector[\Sapprox_i^\perp](\sum_{j= 1}^{\dimSpspcapSop} c_j \orthSoSpsbasis_j + \sum_{j=1}^\dimSo a_j \Sobasis_j^*)},\nonumber\\
&\leq \kvvbar{\sum_{j= 1}^{\dimSpspcapSop} c_j \orthSoSpsbasis_j + \sum_{j=1}^\dimSo a_j \Sobasis_j^*}, \nonumber\\
&= \kparen{\sum_{j= 1}^{\dimSpspcapSop} c_j^2 + \sum_{j=1}^\dimSo a_j^2}^{\frac{1}{2}},\nonumber
\end{align}
where 
the last equality results from the orthogonality of $\kbrace{\Sobasis_i^*}_{i=1}^\dimSo$ and $\kbrace{\orthSoSpsbasis_i}_{i=1}^{\dimSpspcapSop}$. 

 Since the coefficients $a_i$ and $c_i$ must satisfy the constraints in \eqref{eq:defmanifoldperp2}, we finally obtain
\begin{align}
\max_{\vech\in \Mposts} \dist[\vech]{\Sapprox_i} 
&\leq \kparen{\Spswidth^2 + \Spwidth^2}^{\frac{1}{2}}
\leq \Spswidth + \Spwidth. \label{eq:bound3}
\end{align}

We note that since $\svdecr_j=1$ for $j=1,\ldots, \maxosv$, the upper bounds stated in \eqref{eq:bound2} and \eqref{eq:bound3} can be jointly rewritten as
\begin{align}
\max_{\vech\in \Mposts} \dist[\vech]{\Sapprox_i} &\leq \svdecr_{\maxzsv-(i-\dimSp^*)}^{-1} \kparen{\Spswidth+\Spwidth},\label{eq:boundjoint}
\end{align}
for $i =\dimSp^*, \ldots, \dimSp^*+\maxzsv-1$.

\subsubsection{Definition of $\Kspecbound_i$ and $\Kspecboundb_i$}:
$\Kspecboundb_i$ is a direct consequence of the bound derived in Section \ref{sec:defKbarbar}. The definition of $\Kspecbound_i$ results from a combination of the bounds \eqref{eq:bound1} and \eqref{eq:boundjoint} exposed in Sections \ref{sec:defKbar1} and  \ref{sec:defKbar3}.

In order to obtain the tightest bound on $\Kspec_i(\Mpost)$,  $\Kspecbound_i$ must be equal to the minimum of \eqref{eq:bound1} and \eqref{eq:boundjoint} for each index $i$ at which both bounds are defined.
 Since $\Spwidth+\Spswidth\geq \Spswidth$, the bound in \eqref{eq:bound1} is always smaller than \eqref{eq:boundjoint} when $i\geq \dimSps$. This leads to the last line in \eqref{eq:UpperBoundKW1}. On the other hand, when $i< \dimSps$ the bound in \eqref{eq:bound1} is infinite and \eqref{eq:boundjoint} thus takes the lead. This results in the second line in \eqref{eq:UpperBoundKW1}. 
 
 We note that the bound in \eqref{eq:boundjoint} is always well-defined  in the range $i=\dimSp^*,\ldots, \dimSps-1$ since
\begin{align}
\dimSp^*+\maxzsv-1
&= \min\kparen{\dimSps+\maxzsv-1,\dimSps+\dimSp-1}\nonumber,\\
&\geq \dimSps.\nonumber
\end{align}

\section{Conclusions}

In this paper, we tackle the problem of finding a good approximation subspace for a solution manifold $\Ms$. Unlike in the standard setup where the solution manifold is assumed to be known, we assume that only partial information is available on the latter. More specifically, we suppose that 
 we have the following information at our disposal: \textit{i)} we know that the target manifold is included in a larger set, dubbed ``prior manifold''; \textit{ii)} we have access to a set of partial linear observations for each element of $\Ms$. This setup corresponds, for example, to the ubiquitous situation where some parameters of the system to approximate are imperfectly known but some sensing device can provide us with partial measurements of the state of the system. 

In this work, we thus address the following questions: how to combine the prior knowledge and the collected measurements to build a good approximation subspace for $\Ms$? In particular, what performance can one expect? 

We provide an answer to these questions at both a practical and theoretical level.~From a worst-case perspective, we show that the best-achievable performance is characterized by the Kolmogorov width of a well-defined manifold, the so-called ``posterior manifold''.~Motivated by this finding, we propose a tractable algorithm, combining samples from the posterior manifold and a greedy procedure, to achieve performance close to the optimal solution.~The theoretical behavior of the proposed methodology is finally studied in a simplified scenario, where the prior manifold is assumed to be a degenerate ellipsoid.~We emphasize that the performance achievable in the partially-informed setup is highly dependent on the behavior of the singular values of the projector between the prior and the observation subspaces.


\appendix
\section{Technical Details} \label{app:caractMpost}

\subsection{Space Decomposition in the Suitable Bases}\label{sec:charac_MpcapHplan}

In this section we elaborate on some simple facts that will be useful in the rest of the appendix. The material presented hereafter takes the form of two remarks (see Remarks A and B below) and mainly derives from the discussion carried out in \cite[section 2.4]{Binev2015Data}. In particular, we exploit the following lemma from \cite{Binev2015Data}:
\begin{lemma} \label{lemma:ortho_decomposition}
Let $X$ and $Y$ be two linear subspaces of $\Sh$. Then, 
\begin{align}
Y = \projector[Y](X) \oplus \kparen{Y\cap X^\perp}, \label{eq:ortho_decomposition}
\end{align}
where $\projector[Y](X)$ and $\kparen{Y\cap X^\perp}$ are orthogonal. 
\end{lemma}
%
\textbf{Remark A:} 
$\Sh$ can be decomposed as the direct sum of four orthogonal subspaces, namely:
\begin{align}
\Sh 
&= \underbrace{\projector[\So](\Sps) \oplus \kparen{\So\cap \Spsp}}_{=\So}
 \oplus \underbrace{\projector[\Sop](\Sps) \oplus \kparen{\Sop\cap \Spsp}}_{=\Sop}. \label{eq:orthdecompositionH}
\end{align}
This follows from the application of Lemma \ref{lemma:ortho_decomposition} with $Y\leftarrow \So$, $X\leftarrow \Sps$ and $Y\leftarrow \Sop$, $X\leftarrow \Sps$ respectively. 

Moreover, ONBs of the first three spaces in \eqref{eq:orthdecompositionH} can be derived as a function of the suitable bases $\kbrace{\Sobasis_j^*}_{j=1}^\dimSo$, $\kbrace{\Spsbasis_j^*}_{j=1}^\dimSps$ defined in Section \ref{sec:favorable_bases}. First, $\kbrace{\Sobasis^*_j}_{j=1}^{\maxzsv}$  is an ONB of $\projector[\So](\Sps)$ by definition.\footnote{We remind the reader that $\maxzsv$ is defined as $\maxzsv\triangleq \card[ \kbrace{j :\svdecr_j>0}]$.} Second, since $\kbrace{\Sobasis^*_j}_{j=1}^{\dimSo}$ is an ONB for $\So$ and, from Lemma \ref{lemma:ortho_decomposition}, $\So=\projector[\So](\Sps) \oplus \kparen{\So\cap \Spsp}$ with $\projector[\So](\Sps)$ orthogonal to $\So\cap \Spsp$, we have that  $\kbrace{\Sobasis^*_j}_{j=\maxzsv+1}^{\dimSo}$ is an ONB of $\So \cap \Spsp$. 
 Third, 
\begin{align}
\kbrace{\kparen{1-\svdecr_j^2}^{-\frac{1}{2}}\kparen{\Spsbasis^*_j - \svdecr_j \Sobasis^*_j}}_{j=\maxosv+1}^\maxzsv \cup \kbrace{\Spsbasis^*_j}_{j=\maxzsv+1}^\dimSps,
\end{align}
is an ONB of $\projector[\Sop](\Sps)$.
 This can be found by projecting $\kbrace{\Spsbasis_j^*}_{j=1}^\dimSps$ onto $\Sop$ and observing that the resulting elements are mutually orthogonal. More specifically, exploiting property \eqref{eq:correlationoptbases} and using the fact that $\Sobasis_j^*=\Spsbasis_j^*$ when $\svdecr_j=1$, we obtain:
\begin{align}
\projector[\Sop]\kparen{\Spsbasis_j^*} 
&= \Spsbasis_j^* - \projector[\So]\kparen{\Spsbasis_j^*},\nonumber\\
&= 
\left\{
\begin{array}{ll}
\boldsymbol{0} & \mbox{if $j = 1, \ldots, \maxosv$},\\
\Spsbasis^*_j - \svdecr_j \Sobasis^*_j & \mbox{if $j=\maxosv+1, \ldots, \maxzsv$}, \\
\Spsbasis^*_j & \mbox{if $j=\maxzsv+1, \ldots, \dimSps$}.
\end{array}
\right. \label{eq:derivBasisPSop(Sps)}
\end{align}
Let us note that there is usually no ONB of $\Sop \cap \Spsp$ which can be expressed as a simple function of $\kbrace{\Sobasis_j^*}_{j=1}^\dimSo$ and $\kbrace{\Spsbasis_j^*}_{j=1}^\dimSps$. In the sequel, we will thus use the notation $\kbrace{\orthSoSpsbasis_j}_{j=1}^{\dimSpspcapSop}$ to denote an arbitrary ONB of $\Sop \cap \Spsp$.\\[0.2cm]

\textbf{Remark B:}  Using Lemma \ref{lemma:ortho_decomposition}, the subspace $\Spsp$ can be decomposed as the following direct sum: 
\begin{align}
\Spsp = \projector[\Spsp]\kparen{\So} \oplus  \kparen{\Sop \cap \Spsp},
\end{align}
where $\projector[\Spsp]\kparen{\So}$ and $\kparen{\Sop \cap \Spsp}$ are orthogonal subspaces. 
We also have that 
\begin{align}
\kbrace{\kparen{1-\svdecr_j^2}^{-\frac{1}{2}}\kparen{\Sobasis^*_j - \svdecr_j \Spsbasis^*_j}}_{j=\maxosv+1}^\maxzsv \cup \kbrace{\Sobasis^*_j}_{j=\maxzsv+1}^\dimSo,
\end{align}
 constitutes an ONB of $\projector[\Spsp]\kparen{\So}$. This can be seen by projecting the ONB $\kbrace{\Sobasis_j^*}_{j=1}^\dimSo$ onto $\Spsp$ and noticing that the resulting  vectors are mutually orthogonal. 
  In particular, exploiting property \eqref{eq:correlationoptbases} and using the fact that $\Sobasis_j^*=\Spsbasis_j^*$ when $\svdecr_j=1$, we find:
\begin{align}
\projector[\Spsp](\Sobasis^*_j) 
&= \Sobasis_j^* - \projector[\Sps]\kparen{\Sobasis_j^*},\nonumber\\
&= \label{eq:PWontoSpsp}
\left\{
\begin{array}{ll}
\0 & \mbox{for $j=1,\ldots, \maxosv$},\\
\Sobasis^*_j-\svdecr_i \Spsbasis^*_j & \mbox{for $j=\maxosv+1,\ldots, \maxzsv$},\\
\Sobasis_j^* & \mbox{for $j=\maxzsv+1, \ldots,\dimSo$}.
\end{array}
\right.
\end{align}


\subsection{Characterization of $\Mps\cap\Hplan_\vech$}\label{app:def_of_MpscapHplan}

In this section, we provide a mathematical characterization of $\Mps\cap\Hplan_\vech$ in the case where $\Mps$ is defined as in \eqref{eq:PriorSS}. 

Let $\vech\in \Sh$. By definition, for any $\vech'\in\Sh$, we have
\begin{align}
\begin{array}{lcl}
\vech'\in\Hplan_\vech & \iff & \scap[\Sobasis^*_i]{\vech'}=\scap[\Sobasis^*_i]{\vech} \mbox{ for all $i\leq m$},\\
\vech'\in\Mps & \iff & \kvvbar{\projector[\Spsp](\vech')}^2\leq \Spswidth^2. 
\end{array}
\label{eq:constraintHMps}
\end{align}
Because of our considerations in Remark A in Section \ref{sec:charac_MpcapHplan}, any $\vech'\in\Sh$ can be expressed as (in the expression below we simply express $\vech'\in\Sh$ in the ONBs of four subspaces appearing in \eqref{eq:orthdecompositionH}) 
\begin{align}
\vech' 
&= \sum_{j=1}^\dimSo \alpha_j\, \Sobasis_j^* + \sum_{j=\maxosv+1}^\maxzsv \beta_j\, \kparen{1-\svdecr_j^2}^{-\frac{1}{2}}\kparen{\Spsbasis^*_j - \svdecr_j \Sobasis^*_j}
 +\sum_{j=1}^{\dimSpspcapSop} c_j \orthSoSpsbasis_j+\sum_{j=\maxzsv+1}^\dimSps d_j \Spsbasis_{j}^*, \label{eq:decomposition_vech}
\end{align}
for some $\alpha_j$, $\beta_j$, $c_j$, $d_j\in\Rbb$. Hereafter, we express the conditions in \eqref{eq:constraintHMps} in terms of these parameters. 

 On the one hand, the first condition in \eqref{eq:constraintHMps} simply reads
\begin{align}
\alpha_j &= \scap[\Sobasis^*_j]{\vech} \mbox{ for all $j\leq m$}. \label{eq:constraintMpsHplan1}
\end{align}
On the other hand, the second condition in \eqref{eq:constraintHMps} can be rewritten as follows.~Starting from \eqref{eq:decomposition_vech} and 
 using the fact that
\begin{align}\nonumber
\projector[\Spsp](\vech') = \vech' - \sum_{i=1}^{\dimSps} \scap[\Spsbasis^*_j]{\vech'} \Spsbasis_j^*,
\end{align}
we find: 
\begin{align}
\projector[\Spsp](\vech') 
&= \sum_{j=\maxosv+1}^{\maxzsv} \alpha_j \kparen{\Sobasis^*_j-\svdecr_j \Spsbasis^*_j}+ \sum_{j=\maxzsv+1}^\dimSo \alpha_j \Sobasis_j^*
- \sum_{j=\maxosv+1}^{\maxzsv} \beta_j \svdecr_j \kparen{1-\svdecr_j^2}^{-\frac{1}{2}} \kparen{\Sobasis^*_j-\svdecr_j \Spsbasis^*_j}\nonumber\\
&+ \sum_{j= 1}^{\dimSpspcapSop} c_j \orthSoSpsbasis_j, \label{eq:proj_vech_Spsp}\\
&= \sum_{j=\maxosv+1}^{\maxzsv} \kparen{\alpha_j-\beta_j \svdecr_j\kparen{1-\svdecr_j^2}^{-\frac{1}{2}}} \kparen{\Sobasis^*_j-\svdecr_j \Spsbasis^*_j}
+ \sum_{j=\maxzsv+1}^\dimSo \alpha_j \Sobasis_j^* +\sum_{j= 1}^{\dimSpspcapSop} c_j \orthSoSpsbasis_j. \label{eq:proj_vech_Spsp2}
\end{align}
Now, because all the terms in the last expression are orthogonal\footnote[1]{This can easily be seen from Remark B in section \ref{sec:charac_MpcapHplan}.}, the second constraint in \eqref{eq:constraintHMps} takes the form:
 \begin{align}
\Spswidth^2&\geq 
 \sum_{j=\maxosv+1}^{\maxzsv} \kparen{\alpha_j -\beta_j \svdecr_j \kparen{1-\svdecr_j^2}^{-\frac{1}{2}} }^2 \kparen{1-\svdecr_j^2}
+ \sum_{j=\maxzsv+1}^\dimSo \alpha_j^2 + \sum_{j= 1}^\dimSpspcapSop c_j^2,\nonumber\\
&= \sum_{j=\maxosv+1}^{\maxzsv} \svdecr_j^{2} \kparen{\alpha_j \svdecr_j^{-1} \kparen{1-\svdecr_j^2}^{\frac{1}{2}}-\beta_j }^2 
+ \sum_{j=\maxzsv+1}^\dimSo \alpha_j^2 + \sum_{j= 1}^\dimSpspcapSop c_j^2. \label{eq:norm_proj_vech_Spsp}
\end{align}

As a conclusion, \eqref{eq:decomposition_vech} together with \eqref{eq:constraintMpsHplan1} and \eqref{eq:norm_proj_vech_Spsp} fully specify the elements of $\Mps\cap\Hplan_\vech$. 
 For future use, we re-express this system of equations by making the following change of variable: $b_j = \svdecr_j \kparen{\alpha_j \svdecr_j^{-1} \kparen{1-\svdecr_j^2}^{\frac{1}{2}}-\beta_j}$ for $j=\maxosv+1,\ldots,\maxzsv$. We then obtain that any $\vech'\in\Mps\cap\Hplan_\vech$ can be written as\footnote{To obtain the first term in \eqref{eq:decomposition_vech2}, we use the fact that $\sum_{j=\maxosv+1}^\maxzsv \alpha_j \svdecr_j^{-1} \Spsbasis^*_j + \sum_{j=1}^\maxosv \alpha_j \Sobasis^*_j = \sum_{j=1}^\maxzsv \alpha_j \svdecr_j^{-1} \Spsbasis^*_j$ since $\Spsbasis^*_j = \Sobasis^*_j$ and $\svdecr_j=1$ for $j\leq \maxosv$, see Section \ref{sec:favorable_bases}.}
\begin{align}
\vech' 
&= \sum_{j=1}^\maxzsv \alpha_j \svdecr_j^{-1} \Spsbasis^*_j + \sum_{j=\maxzsv+1}^\dimSo \alpha_j \Sobasis_j^*
- \sum_{j=\maxosv+1}^{\maxzsv} b_j \svdecr_j^{-1}\kparen{1-\svdecr_j^2}^{-\frac{1}{2}}\kparen{\Spsbasis^*_j - \svdecr_j \Sobasis^*_j}\nonumber\\
&+ \sum_{j= 1}^{\dimSpspcapSop} c_j \orthSoSpsbasis_j + \sum_{j=\maxzsv+1}^{\dimSps} d_j \Spsbasis^*_j, \label{eq:decomposition_vech2}
\end{align}
with 
\begin{align} \label{eq:constrainttot}
\left\{
\begin{array}{l}
\alpha_j = \scap[\Sobasis^*_j]{\vech} \mbox{ for all $j\leq m$},\\
\sum_{j=\maxzsv+1}^\dimSo \alpha_j^2 +\sum_{j= \maxosv+1}^\maxzsv b_j^2 +\sum_{j= 1}^{\dimSpspcapSop} c_j^2 \leq \Spswidth^2.
\end{array}
\right.
\end{align}

\subsection{Characterization of $\cup_{\vech\in\Mp} \kparen{\Mps\cap \Hplan_\vech}$}

In this section we show that $\cup_{\vech\in\Mp} \kparen{\Mps\cap \Hplan_\vech}$ can be expressed as stated in \eqref{eq:formunionEllipsoid}-\eqref{eq:defmanifoldperp2}.~We start from \eqref{eq:decomposition_vech2}-\eqref{eq:constrainttot} and particularize the expression of the $\alpha_i$'s when $\vech\in\Mp$ to obtain the result. 

First, notice that any $\vech\in\Mp$ can (by definition) be written as
$\vech = \uu + \zz$
where $\uu\in\Sp$, $\zz\in B_{\Spwidth}$ and $B_{\Spwidth}=\kbrace{\vech\in\Sh: \kvvbar{\vech}\leq\Spwidth}$ is the $\kvvbar{\cdot}$-ball of radius $\Spwidth$.   
Using the fact that $\Sp\subseteq\Sps$ and $\Sh = \So \oplus \Sop$,  $\uu$ and $\zz$ can be re-expressed as 
\begin{align}
\uu &= \sum_{j=1}^\dimSps \scap[\Spsbasis_j^*]{\uu} \Spsbasis^*_j,\\
\zz &= \sum_{j=1}^\dimSo a_j \Sobasis^*_j + \projector[\Sop](\zz),
\end{align}
we obtain from \eqref{eq:correlationoptbases} that
\begin{align}
\alpha_j &\triangleq \scap[\Sobasis_j^*]{\vech},\nonumber\\
&= \scap[\Sobasis_j^*]{\uu}+\scap[\Sobasis_j^*]{\zz},\nonumber\\
&= \left\{
\begin{array}{ll}
\svdecr_j \scap[\Spsbasis_j^*]{\uu}+a_j & \mbox{if $j=1,\ldots,\maxzsv$,} \\
a_j & \mbox{if $j=\maxzsv+1, \ldots, \dimSo$.}
\end{array} 
\right.
\end{align}
Hence, the first term in \eqref{eq:decomposition_vech2} becomes:
\begin{align}
\sum_{j=1}^\maxzsv \alpha_j \svdecr_j^{-1} \Spsbasis^*_j 
&= \sum_{j=1}^\maxzsv \scap[\Spsbasis_j^*]{\uu} \Spsbasis^*_j + \sum_{j=1}^\maxzsv  \svdecr_j^{-1} a_j \Spsbasis^*_j.
\end{align}
Moreover, because $\zz\in B_{\Spwidth}$, the $a_j$'s must verify:
\begin{align}\nonumber
\sum_{j=1}^\dimSo a_j^2 \leq \Spwidth^2, 
\end{align}
since
\begin{align}\nonumber
\sum_{j=1}^\dimSo a_j^2 = \kvvbar{\projector[\So](\zz)}^2\leq \kvvbar{\zz}^2 \leq \Spwidth^2. 
\end{align}

Combining these results and re-arranging the terms in \eqref{eq:decomposition_vech2}, we obtain that $\cup_{\vech\in\Mp} \kparen{\Mps\cap \Hplan_\vech}$ is characterized by the following set of equations:
\begin{align}
\vech' &
= \sum_{j=1}^\maxzsv \scap[\Spsbasis_j^*]{\uu} \Spsbasis^*_j 
+ \sum_{j=\maxzsv+1}^{\dimSps} d_j \Spsbasis^*_j
- \sum_{j=\maxosv+1}^{\maxzsv} b_j \svdecr_j^{-1}\kparen{1-\svdecr_j^2}^{-\frac{1}{2}}\kparen{\Spsbasis^*_j - \svdecr_j \Sobasis^*_j}\nonumber\\
& + \sum_{j= 1}^{\dimSpspcapSop} c_j \orthSoSpsbasis_j + \sum_{j=\maxzsv+1}^\dimSo a_j \Sobasis_j^*+ \sum_{j=1}^\maxzsv a_j \svdecr_j^{-1} \Spsbasis^*_j, \label{eq:decomposition_vech3}
\end{align}
with
\begin{align} \label{eq:constrainttot2}
\left\{
\begin{array}{l}
\sum_{j=\maxzsv+1}^\dimSo a_j^2 +\sum_{j= \maxosv+1}^\maxzsv b_j^2 +\sum_{j= 1}^{\dimSpspcapSop} c_j^2 \leq \Spswidth^2,\\
\sum_{j=1}^\dimSo a_j^2 \leq \Spwidth^2. 
\end{array}
\right.
\end{align}
Because $\uu\in\Sp$ and there is no constraint on the $d_j$'s, the first two terms in \eqref{eq:decomposition_vech3} define the subspace $\projector[\Sps](\Sp)\oplus \spa[\kbrace{\Spsbasis_j^*}_{j=\maxzsv+1}^\dimSps]=\Sp \oplus \spa[\kbrace{\Spsbasis_j^*}_{j=\maxzsv+1}^\dimSps]$.

It thus remains to show that the last terms in \eqref{eq:decomposition_vech3} together with the constraints \eqref{eq:constrainttot2} define $\manifoldperp$ as in \eqref{eq:defmanifoldperp1}-\eqref{eq:defmanifoldperp2}. This simply follows by re-expressing the last term in \eqref{eq:decomposition_vech3}, that is $\sum_{j=1}^\maxzsv a_j \svdecr_j^{-1} \Spsbasis^*_j$, in terms of its components in $\So$ and $\Sop$. More specifically, using the suitable bases for $\projector[\So](\Sps)$ and $\projector[\Sop](\Sps)$ defined in Section \ref{sec:charac_MpcapHplan}, we find 
\begin{align}
\projector[\So](\sum_{j=1}^\maxzsv  \svdecr_j^{-1} a_j \Spsbasis^*_j)&=\sum_{j=1}^\maxzsv  a_j \Sobasis^*_j,\nonumber\\
\projector[\Sop](\sum_{j=1}^\maxzsv  \svdecr_j^{-1} a_j \Spsbasis^*_j)&= \sum_{j=\maxosv+1}^\maxzsv  a_j \svdecr_j^{-1} \kparen{\Spsbasis^*_j-\svdecr_j \Sobasis^*_j}.\nonumber
\end{align}
The result then follows by plugging these expressions into \eqref{eq:decomposition_vech3}. 

\section{Weak and algebraic formulations of the Thermal-block Problem}\label{app:wformualtion}

Let $H^1(\Omega)$ denotes the Sobolev space of order 1 on $\Omega$ and let
\begin{align}
X = \kbrace{\vech \in H^1(\Omega): \mbox{$\vech$ satisfies the Dirichlet conditions of problem \eqref{eq:thermalBlockproblem}}}. \nonumber
\end{align}
Using the derivations of \cite[Section 2.3.1]{Quarteroni2016Reduced}, we have that the weak formulation of \eqref{eq:thermalBlockproblem} can be written as
\begin{align}
\mbox{find $\vech\in X$ such that\ $a(\vech,\vech') = b(\vech')\ \forall \vech'\in X$}\label{eq:weakformulation1}
\end{align}
where 
\begin{align}
a(\vech,\vech') &= \int_{\Omega} k(\x,\param) \ktranspose{\nabla} \vech  \nabla \vech' d\x,\nonumber\\
b(\vech')          &= \int_{\Omega} s(\x) \vech' d\x+ \int_{\Gamma_1} c\, k(\x,\param)^{-1} \vech' d\x.  \nonumber
\end{align}

Let $X_{\mathrm{fe}}$ be a finite-dimensional subspace of $X$ (for example generated by a basis of finite elements). A discrete approximation of \eqref{eq:weakformulation1} can then be defined as
\begin{align}
\mbox{find $\vech \in X_{\mathrm{fe}}$ such that\ $a(\vech,\vech') = b(\vech')\ \forall \vech'\in X_{\mathrm{fe}}$}.\label{eq:weakformulationapprox}
\end{align}
Considering an ONB of $X_{\mathrm{fe}}$, say $\kbrace{\varphi_j}_j^{\mathrm{N}_\mathrm{fe}}$, the algebraic formulation of \eqref{eq:weakformulationapprox} reads
\begin{align}
\mbox{find the solution of $\A \h = \mathbf{b}$},\label{eq:weakformulationalgebraic}
\end{align}
where the elements of matrix $\A\in\Rbb^{\mathrm{N}_\mathrm{fe} \times \mathrm{N}_\mathrm{fe}}$ and vector $\mathbf{b}\in\Rbb^{\mathrm{N}_\mathrm{fe}}$ are defined as
\begin{align}
A_{i,j} &= a(\varphi_i,\varphi_j),\nonumber\\
b_{j} &= b(\varphi_j).\nonumber
\end{align}
We note in particular that if $c=0$ and the source term decomposes as 
\begin{align}
s(\x) = \sum_j s_j \varphi_j, \label{eq:decsource}
\end{align}
we simply have $b_j = s_j$ because of the orthogonality of the $\varphi_j's$.

Let us now elaborate briefly on the second simulation setup considered in Section \ref{sec:simuresults}. We argue that any vector $\h\in\Rbb^{\mathrm{N}_\mathrm{fe}}$ (or equivalently any $\vech\in X_{\mathrm{fe}}$) can be obtained as the solution of \eqref{eq:weakformulationalgebraic} (resp. of \eqref{eq:weakformulationapprox}) by properly choosing the source term $\sourceterm(\x)$. This can be seen as follows. First note that $\A$ is a fixed matrix if the parameter $\param$ is set to a fixed value. Let $\s$ be the vector gathering the $s_j$'s in \eqref{eq:decsource}. Then, by choosing $\mathbf{s} = \A^{-1}\h$,
it can be seen that the solution of \eqref{eq:weakformulationalgebraic} is $\h$ because we assume that $c=0$ in our simulation setup. The solution of \eqref{eq:weakformulationapprox} is obtained by using \eqref{eq:decsource} with $\mathbf{s} = \A^{-1}\h$.

%
%
%
%
%

\section{Summary Appendix}\label{app:summary}

This appendix is intended to ease the reading of the paper by gathering the main elements appearing in the core of the paper. 
The goal is to allow the reader to access ``at a glance'' to the most important definitions and relationships.  
 
The paper addresses the problem of finding a good approximation subspace for a solution manifold $\Ms$ defined as
\begin{align}\nonumber
\Ms = \{\vech(\param) \in \Sh: \param\in\paramSet\},
\end{align}
where $\vech(\param)$ is the solution of some PPDE depending on parameter $\theta$. No particular constraints are imposed on $\Ms$. The best possible worst-case approximation error of the elements $\Ms$ in a subspace of dimension $i$ is given by the Kolmogorov $i$-width:
\begin{align}\nonumber
\Kspec_i(\Ms) &= \inf_{\Sapprox:\dim(\Sapprox)=i} \kparen{\sup_{\vech\in\Ms} \dist[\vech]{S}}. 
\end{align}
If the above optimization problem admits a minimizer, we denote the latter by $\Sperf_i$, \ie
\begin{align}\nonumber
\Sperf_i &\in \kargmin_{\Sapprox:\dim(\Sapprox)=i} \kparen{\sup_{\vech\in\Ms} \dist[\vech]{S}}.
\end{align}

The construction of the proposed approximation subspace is based on two ingredients that we remind here:
\begin{itemize}
\item[\textit{i)}] \textit{a prior manifold $\Mps$}:  the only constraint we impose on $\Mps$ is to be such that
\begin{align}\nonumber
\Ms\subseteq \Mps. 
\end{align}
In the main body of the paper, we repeatedly use two particular instances of $\Mps$: 
\begin{itemize}
\item An intersection of degenerate ellipsoids, \ie
\begin{align}\nonumber
\Mps=\cap_{j=1}^{\nsubspace} \kbrace{\vech : \dist[\vech]{\Sps_j}\leq \Spswidth_j}.
\end{align}
\item A single degenerate ellipsoid, \ie
\begin{align}\nonumber
\Mps=\kbrace{\vech : \dist[\vech]{\Sps}\leq \Spswidth}.
\end{align}
\end{itemize}

\item[\textit{ii)}]  \textit{a set of partial observations of $\Ms$}: we assume that we collect, $\forall\vech\in\Ms$,  a set of noiseless linear measurements:
 \begin{align} \nonumber
\kbrace{\scap[\Sobasis_j]{\vech} }_{j=1}^\dimSo,
\end{align}
 for some orthonormal basis  $\kbrace{\Sobasis_j }_{j=1}^\dimSo$. We denote by  $\So$ the  $\dimSo$-dimensional subspace induced by $\kbrace{\Sobasis_j }_{j=1}^\dimSo$, \ie
 \begin{align}\nonumber
\So = \spa[\kbrace{\Sobasis_j }_{j=1}^\dimSo].
\end{align}
Moreover, for a given $\vech\in\Ms$, we denote by $\Hplan_\vech$ the set of element leading to the same measurement sets, \ie
\begin{align}
\Hplan_\vech 
&= \kbrace{\vech' : \langle \Sobasis_j, \vech' \rangle =\langle \Sobasis_j, \vech \rangle \mbox{ for $j=1,\ldots, m$}}, \nonumber\\
&= \kbrace{\vech'= \vech +  \Sobasis^\perp: \Sobasis^\perp\in \So^\perp}. \nonumber
\end{align}
\end{itemize}

Our procedure to find a good approximation subspace from the above ingredients is based on the construction of a ``posterior'' manifold whose definition is given by
\begin{align}\nonumber
\Mpost = \Mps \cap \kparen{\cup_{\vech\in\Ms} \Hplan_\vech}. 
\end{align}
By definition,  $\Mpost$ obeys the following relationship
\begin{align}\nonumber
\Ms \subseteq \Mpost \subseteq \Mps. 
\end{align}
We denote by $\Spost_i$ the $i$-dimensional subspace (assuming it exists) minimizing the worst approximation error over the elements of $\Mpost$, \ie
\begin{align}\nonumber
\Spost_i &\in \kargmin_{\Sapprox: \dim(\Sapprox)=i}\kparen{\sup_{\vech\in\Mpost} \dist[\vech]{\Sapprox}}.
\end{align}
This subspace verifies the following property:
\begin{align}\nonumber
\Kspec_i(\Ms) \leq \sup_{\vech\in\Ms} \dist[\vech]{\Spost_i} \leq \Kspec_i(\Mpost)\leq \Kspec_i(\Mps). 
\end{align}

Finally, in many parts of the paper, we exploit the following decomposition of $\Sh$ in four orthogonal subspaces:
\begin{align}\nonumber
\Sh &= \underbrace{\projector[\So](\Sps) \oplus \kparen{\So\cap \Spsp}}_{=\So}
 \oplus \underbrace{\projector[\Sop](\Sps) \oplus \kparen{\Sop\cap \Spsp}}_{=\Sop}.
\end{align}
We also use the following  ONBs for each subspace appearing in the above decomposition:
\begin{align}\nonumber
\begin{array}{rll}
\bullet &\mbox{ONB for $\projector[\So](\Sps)$:} & \kbrace{\Sobasis^*_j}_{j=1}^{\maxzsv},\\
\bullet &\mbox{ONB for $\kparen{\So\cap \Spsp}$:} & \kbrace{\Sobasis^*_j}_{j=\maxzsv+1}^{\dimSo},\\
\bullet &\mbox{ONB for $\projector[\Sop](\Sps)$:} &  \kbrace{\kparen{1-\svdecr_j^2}^{-\frac{1}{2}}\kparen{\Spsbasis^*_j - \svdecr_j \Sobasis^*_j}}_{j=\maxosv+1}^\maxzsv \cup \kbrace{\Spsbasis^*_j}_{j=\maxzsv+1}^\dimSps,\\
\bullet &\mbox{ONB for $\kparen{\Sop\cap \Spsp}$:} & \kbrace{\orthSoSpsbasis_j}_{j=1}^{\dimSpspcapSop},
\end{array}
\end{align}
where $\kbrace{\Spsbasis_j^*}_{j=1}^\dimSps$  and $\kbrace{\Sobasis_j^*}_{j=1}^\dimSo$ denotes the ``suitable'' ONBs for $\Sps$ and $\So$ introduced in Section~\ref{sec:favorable_bases} and $\kbrace{\svdecr_j}_{j=1}^\mindim$ are the singular values of the Gram matrix $\G$ defined as $g_{ij} \triangleq \scap[\Sobasis_i]{\Spsbasis_j}$.
 The quantities $p$ and $q$ are defined as
\begin{align}
\begin{array}{ll}
\maxosv & \triangleq  \card[ \kbrace{j :\svdecr_j=1}],\nonumber\\
\maxzsv & \triangleq  \card[ \kbrace{j :\svdecr_j>0}].\nonumber\\
\end{array}
\end{align}




\section{Acknowledgements}
The authors wish to thank Stefano Pagani for his kind advice and his help in using the Matlab\textsuperscript{\textregistered} toolbox ``redbKIT''. 
The authors also thank the ``Agence Nationale de la Recherche" (ANR) which partially funded this research through the GERONIMO project.

\bibliographystyle{IEEEbib}
\bibliography{group-15332,cherzet}

%
%
%
\end{document}

%% file: headings_herzet.tex
\usepackage[T1]{fontenc} 
\usepackage[latin1]{inputenc}
\usepackage[frenchb,english]{babel}

\usepackage{amssymb} 
\usepackage{amsmath} 
\usepackage{amsfonts} 
\usepackage{mathrsfs}
\usepackage{amsthm}
\usepackage{kmath} 			
\usepackage[noadjust]{cite}       
\usepackage{soul}		         

\usepackage{graphicx}
\usepackage{color}
\usepackage{pgfplots}
\usepackage{epsfig}
\usepackage{epstopdf}

\def\c{{\mathbf c}}

\def\h{{\mathbf h}}

\def\n{{\mathbf n}}

\def\s{{\mathbf s}}

\def\x{{\mathbf x}}

\def\0{{\mathbf 0}}
\def\1{{\mathbf 1}}

\def\A{{\mathbf A}}

\def\G{{\mathbf G}}

\def\S{{\mathbf S}}

\def\X{{\mathbf X}}

\def\Z{{\mathbf Z}}

\def\hh{{\boldsymbol h}}

\def\uu{{\boldsymbol u}}
\def\vv{{\boldsymbol v}}
\def\ww{{\boldsymbol w}}

\def\zz{{\boldsymbol z}}

\def\Ec{{\mathcal E}}

\def\Hc{{\mathcal H}}

\def\Mc{{\mathcal M}}

\def\Oc{{\mathcal O}}

\def\Rbb{{\mathbb R}}

\def\ie{\textit{i.e.},\xspace}
\def\eg{\textit{e.g.},\xspace}
\def\etal{\textit{et al.}\xspace }

\newkmacro\sign[1][{\cdot}]{\mathrm{sign}\left({#1}\right)} 
\newkmacro\interior[1][{\cdot}]{\mathrm{int}\left({#1}\right)}
\newkmacro\bd[1][{\cdot}]{\mathrm{bd}\left({#1}\right)}
\newkmacro\conv[1][{\cdot}]{\mathrm{conv}\left({#1}\right)}
\newkmacro\card[1][{\cdot}]{\mathrm{card}\left({#1}\right)}
\newkmacro\spa[1][{\cdot}]{\mathrm{span}\left({#1}\right)}
\newkmacro\nul[1][{\cdot}]{\mathrm{null}\left({#1}\right)}
\newkmacro\deter[1][{\cdot}]{\mathrm{det}\left({#1}\right)}
\newkmacro\tr[1][{\cdot}]{\mathrm{tr}\left({#1}\right)}

\newtheorem{lemma}{Lemma}
\newtheorem{theorem}{Theorem}

\newtheorem{example}{Example}


%% file: def_macro.tex
\def\Sh{{\Hc}} 				
\def\vech{{\hh}} 		         
\def\vechtrial{{\hh'}} 		         
\def\dimSh{{\dim(\Sh)}} 		 

\def\dimSp{{k}} 		        		
\def\Sp{{U}}				
\def\Spwidth{{\epsilon}}		
\def\Spbasis{{\uu}} 		        
\def\Mp{{\Sigma}}		        
\def\Ms{{\Mc}}				
\def\Mstrial{{\tilde{\Ms}}}

\def\Spext{{\hat{\Sp}}}		

\def\dimSps{{n}} 		        	
\def\Sps{{V}}				
\def\Spswidth{{\hat{\Spwidth}}}      
\def\Spsbasis{{\vv}} 		        
\def\Spsp{{\Sps}^\perp}		
\def\Mps{{\Mc_{\mathrm{prior}}}}		 
\def\Spsbasist{{\tilde{\Spsbasis}}} 		        
\def\dimSpsmax{{\mathrm{N}_{\mathrm{max}}}}

\def\dimSo{{m}} 		        
\def\So{{W}} 				
\def\Sop{{\So}^\perp} 		
\def\Sobasis{{\ww}} 		        
\def\Mpost{{\Ms_{\mathrm{post}}}}			
\def\Mposts{{\hat{\Ms}_{\mathrm{post}}}}		
\def\intbasis{{\boldsymbol{\phi}}}				
\def\Sobasist{{\tilde{\Sobasis}}} 		        

\def\Hplan{{H}}				
\def\Kspec{{\kappa}}			
\def\Kspecbound{{\bar{\Kspec}}}        
\def\Kspecboundb{{\bar{\bar{\Kspec}}}}        
\def\orthSoSpsbasis{{\boldsymbol{\psi}}}	
\def\maxosv{{p}}			
\def\maxzsv{{q}}			
\def\mindim{{\min\kparen{\dimSo,\dimSps}}} 

\def\manifoldperp{{\Ec}}
\newkmacro{\dist}[2][]{\mathrm{dist}(#1,#2)}
\newkmacro{\scap}[2][]{\kangle{#1,#2}}
\newkfunc{\projector}{P}

\def\Sapprox{{S}}
\def\Sapproxa{{\hat{\Sapprox}}}
\def\Sperf{\Sapprox^{\mathrm{perf}}}  
\def\Spost{\Sapprox^{\mathrm{post}}} 
\def\Spointa{\hat{\Sapprox}^{\mathrm{point}}}  

\def\param{{\theta}}
\def\paramstep{{\param_\mathrm{step}}}
\def\paramSet{{\Theta}}
\def\paramSetp{{\paramSet_{\mathrm{relax}}}}
\def\Mrs{{\Mc_{\mathrm{relax}}}}		 

\def\svdecr{{\lambda}} 		                
\def\dimSpspcapSop{{\dim\kparen{\Sop\cap\Spsp}}}  

\def\feasmanifoldset{{\Xi_{\mathrm{feas}}}}
\def\nsnapshots{{\mathrm{N}_{\mathrm{snap}}}}
\def\nsamples{{\mathrm{N}_{\mathrm{samples}}}}
\def\nsubspace{{\mathrm{L}}}
\def\vechVpWp{{\zz}}
\def\dstep{{t}}
\def\nstep{{\mathrm{T}}}
\def\cel{{\vech_{\mathrm{center}}}} 				
\def\sourceterm{{s}}
\def\corrVW{{\delta}}